\documentclass[12pt]{article}

\usepackage{amsfonts,epsfig,amssymb,array}

 \newtheorem{thm}{Theorem}
 \newtheorem{propos}[thm]{Proposition}
 \newtheorem{lem}[thm]{Lemma}
 \newtheorem{cor}[thm]{Corollary}

 \newcommand{\eps}{\varepsilon}
 \newcommand{\la}{\lambda}
 
 \newcommand{\ee}{\textup{e}}

 \newcommand{\lep}{\Lambda_{\varepsilon}(P)}
 
 \newcommand{\interior}{\textup{Int}}
 \newcommand{\GG}{\mathcal{G}}
 \newcommand{\PP}{\mathcal{P}}
 \newcommand{\RR}{\mathbb{R}}
 
 \newcommand{\CC}{\mathbb{C}}
 
 \newcommand{\F}{\mathcal{F}}
 \newcommand{\LL}{\mathcal{L}}

 \newcommand{\B}{\mathcal{B}}
 \newcommand{\N}{\mathcal{N}}

\newcommand{\mat}{\left[\begin{array}}
\newcommand{\rix}{\end{array}\right]}
\newcommand{\qed}{\ifmmode$\Box$\else{\unskip\nobreak\hfil
\penalty50\hskip1em\null\nobreak\hfil$\Box$
\parfillskip=0pt\finalhyphendemerits=0\endgraf}\fi}
\newcommand{\eq}{\begin{equation}}
\newcommand{\en}{\end{equation}}
\newcommand{\eqa}{\begin{eqnarray}}
\newcommand{\ena}{\end{eqnarray}}
\begin{document}

\title{On Pseudospectra of Matrix Polynomials \\ and their Boundaries}

\author{Lyonell Boulton, \\
{\small Department of Mathematics and the Maxwell Institute for Mathematical
Sciences,} \\
{\small  Heriot-Watt University, Edinburgh EH14 2AS, United Kingdom.}
    \and
Peter Lancaster, \\
{\small Department of Mathematics and Statistics,} \\
{\small University of Calgary, Calgary AB, Canada T2N 1N4.}
    \and
Panayiotis Psarrakos, \\
{\small Department of Mathematics, National Technical University,} \\
{\small Zografou Campus, 15780 Athens, Greece.}}

\maketitle

% --------------------------------------------------------
\begin{abstract}
In the first part of this paper (Sections~\ref{s2}-\ref{s4}), the
main concern is with the boundary of the pseudospectrum of a
matrix polynomial and, particularly, with smoothness properties of
the boundary. In the second part
(Sections~\ref{number}-\ref{examples}), results are obtained
concerning the number of connected components of pseudospectra,
as well as results concerning matrix polynomials with multiple
eigenvalues, or the proximity to such polynomials.
\end{abstract}

% --------------------------------------------------------

\

\vspace{.1cm}



\pagebreak

%%%%%%%%%%%%%%%%%%%%%%%%%%%%%%%%%%%%%%%%%%%%%%%%%%%%%%%%%%
\section{Introduction} \label{s1}
%%%%%%%%%%%%%%%%%%%%%%%%%%%%%%%%%%%%%%%%%%%%%%%%%%%%%%%%%%

This paper falls into two parts. In the first
(Sections~\ref{s2}-\ref{s4}), the main concern is with the
boundary of the pseudospectrum of a matrix polynomial and,
particularly, in view of its importance for boundary-tracing
algorithms, with the smoothness properties of the boundary. In
the second (Sections~\ref{number}-\ref{examples}), we further
develop analysis begun by two of the present authors (see
\cite{LP}) on qualitative aspects of the pseudospectrum. This
part is also influenced by earlier work on pseudospectra for
standard eigenvalue problems by Alam and Bora in \cite{AB}. In
particular, results are presented concerning the number of
connected components of the pseudospectrum and proximity to
systems with multiple eigenvalues.

Let us begin with some formal definitions. First, a {\it matrix
polynomial} is a function $P:\CC\rightarrow\CC^{n\times n}$ (the
algebra of all $\,n\times n\,$ complex matrices) of the form \eq
\label{eq.poly}
  P(\la) \,=\, P_m \la^m + P_{m-1} \la^{m-1} +
                \cdots + P_1 \la + P_0 ,
\en where $\la$ is a complex variable and $\,P_0,P_1 ,\dots
,P_m\in \CC^{n\times n}\,$ with $\,\det P_m \ne 0$. The {\it
spectrum} of such a function is $\,\sigma(P):= \{\la \in
\CC\,:\,\textup{det}\,P(\la)=0\}$.

Since $\,\det P_m\ne 0,\,$ $\sigma(P)$ consists of no more than
$nm$ distinct \textit{eigenvalues}. A nonzero vector
$\,x_0\in\CC^n\,$ is known as an \textit{eigenvector} of $P(\la)$
corresponding to an eigenvalue $\,\la_0\in\sigma(P)\,$ if it
satisfies $\,P(\la_0)x_0=0$. The \textit{algebraic multiplicity}
of a $\,\la_0\in\sigma(P)\,$ is the multiplicity of $\la_0$ as a
zero of the scalar polynomial $\det P(\la)$, and it is always
greater than or equal to the \textit{geometric multiplicity} of
$\la_0$, that is, the dimension of the null space of the matrix
$P(\la_0)$. A multiple eigenvalue of $P(\la)$ is called
\textit{defective} if its algebraic multiplicity exceeds its
geometric multiplicity.

We let $\PP_m$ denote the linear space of $\,n\times n\,$ matrix
polynomials with degree $m$ or less. Using the spectral matrix
norm (i.e., that norm subordinate to the Euclidean vector norm),
we may define the max norm on $\PP_m$,
\eq  \label{eq.norm}
   \|P(\la)\| \,:=\, \textup{max}_{0\le j\le m} \| P_j \| .
\en Using this norm, we construct a class of matrix polynomials
obtained from $P(\la)$ in (\ref{eq.poly}) by perturbation. The
admissible perturbations are defined in terms of a real polynomial
$\,w(x)=\sum_{j=0}^mw_j x^j\,$ with nonnegative coefficients and a
positive constant coefficient; $\,w_j \ge 0\,$ for each
$\,j=1,2,\dots,m,\,$ and $\,w_0>0$. First consider matrix
polynomials in $\PP_m$ of the form
\begin{equation} \label{e1}
   Q(\la) \,=\, (P_m + \Delta_m) \la^m + \cdots +
                (P_1 + \Delta_1) \la + (P_0 + \Delta_0),
\end{equation}
where the matrices $\,\Delta_j\in\CC^{n\times n}\,$
($j=0,1,\dots,m$) are arbitrary. Then, for a given $\,\eps
\geq 0 ,\,$ the class of admissible perturbed matrix polynomials is
\eq \label{eq.b}
  \B(P,\eps, w) \,:=\, \left \{ Q(\la) :\, \|\Delta_j\|
                      \le \eps\, w_j ,\; j=0,1,\ldots,m \right \} .
\en
This is a convex compact set in the linear space $\PP_m$ with the
norm (\ref{eq.norm}).

The $\eps$-{\em pseudospectrum of} $P(\la)$ {\em with respect to}
$w(x)$ (introduced by Tisseur and Higham \cite{TH}) is then
\eq\label{eq.ps}
  \lep \,:=\, \left \{ \mu \in \CC : \, \det Q(\mu) = 0
    \;\, \mbox{for some}\;\, Q(\la)\in\B(P,\eps,w) \right \}.
\en As $w(x)$ is generally fixed throughout this paper, it
will not appear explicitly in this notation, and we will refer to
$\lep$ simply as the $\eps$-{\em pseudospectrum of} $P(\la)$. Note
that if $\,\eps\, w_m < \|P_m^{-1}\|^{-1},\,$ then all matrix
polynomials in $\B(P,\eps,w)$ have nonsingular leading
coefficients, and this ensures that $\lep$ is bounded (Theorem 2.2
of \cite{LP}).

If we define the {\em standard} eigenvalue problem as that in
which $\,P(\la)=I\la -A,\,$ then it is natural to define weights
$\,w_1=0\,$ (no perturbation of the coefficient $I$ is admitted)
and $\,w_0=1$. Thus, $\,w(x)=1\,$ and, using (\ref{eq.ps}), we
obtain the relatively well-understood ``$\eps$-pseudospectrum of
matrix $A$'', namely,
\[
  \Lambda_{\eps}(A) \,\equiv\, \lep \,=\,
  \left \{\mu \in \CC :\, \det(I \mu -(A+\Delta_0)) = 0,
      \; \|\Delta_0\| \le \eps \right \} .
\]

%%%%%%%%%%%%%%%%%%%%%%%%%%%%%%%%%%%%%%%%%%%%%%%%%%%%%%%%%%
\section{The singular value functions} \label{s2}
%%%%%%%%%%%%%%%%%%%%%%%%%%%%%%%%%%%%%%%%%%%%%%%%%%%%%%%%%%

For any $\,\la\in\CC,\,$ the singular values of a matrix
polynomial $P(\la)$ are the nonnegative square-roots of the $n$
eigenvalue functions of $P(\la)^*P(\la)$. They are denoted by
\[
 s_1(\la)\,\ge\,s_2(\la)\,\ge\,\cdots\,\ge\,s_n(\la)\,\ge\,0.
\]
The real-valued function $\,s_n:\CC\longrightarrow [0,\infty),\,$
given by the smallest singular value, provides more information
about the matrix polynomial $P(\la)$ than $\sigma (P)$ alone.
This will become clear in the forthcoming section when we discuss
the pseudospectrum of $P(\la)$. Let us first describe some
general properties of $s_n(\la)$.

It is clear that an alternate definition of the spectrum of a
matrix polynomial $P(\la)$ is:
\[
  \sigma(P) \,=\, \left \{\la\in\CC :\,s_n(\la)=0 \right \}.
\]
The connection between the zeros of $s_n(\la)$ and the eigenvalues
of $P(\la)$ can be made more precise using the singular value
decomposition.

\begin{propos}\label{prop.a}
An eigenvalue $\,\la_0 \in \sigma(P)\,$ has geometric multiplicity
$\,k$ if and only if
\[
    s_1(\la_0) \,\ge\,  s_2(\la_0) \,\ge\,
    \cdots \,\ge\, s_{n-k}(\la_0) \,>\, s_{n-k+1}(\la_0)
    \,=\, \cdots \,=\, s_n(\la_0) \,=\, 0 .
\]
\end{propos}

Our analysis depends on an important, concise characterisation of
the $\eps$-pseudospectrum in terms of the function $s_n(\la)$.
This was obtained by Tisseur and Higham (Lemma 2.1 of \cite{TH}),
\eq\label{eq.th}
    \lep \,=\, \left \{ \la \in \CC : \,
               s_n(\la) \le \eps \, w(|\la|) \right \} .
\en Clearly, $\,\sigma(P)=\Lambda_0(P) \subset \lep \,$ for any
$\,\eps > 0$. Thus, $\lep$ is nothing but the level set  at height
$0$ of the real-valued function $\,s_n(\la)-\eps\, w(|\la|)\,$, or
that at height $\eps$ of the function $\,s_n(\la)\,
w(|\la|)^{-1}$. Notice also that in the standard eigenvalue
problem, $\,\eps\, w(|\la|)=\eps\,$ in (\ref{eq.th}). More
generally, $\,\eps\,w(|\la|)\,$ (in equation (\ref{eq.th})) is a
radially symmetric non-decreasing function of $\la$.

By using the Euclidean vector norm,
\eq \label{eq.l1}
  s_n(\la) \,=\, \min_{u\ne 0} \frac{\|P(\la)u\|}{\|u\|} \, .
\en Our first theorem has been originally established by Davies in the more general context of holomorphic families of
bounded operators. A proof is included here for completeness.

\begin{thm} \label{the.sa}
Let $P(\la)$ be invertible on a domain $U$. Then $s_n(\la)^{-1}$
is a subharmonic function on $U$.
\end{thm}

\noindent{\bf Proof.} First we recall one of the characterisations
of continuous subharmonic functions (see Ahlfors \cite{A}, for
example). A continuous function $\,\phi :U \rightarrow \RR\,$ is
subharmonic if and only if, for any closed disc in $U$ with centre
$\la_0$ and radius $r$,
%\eq \label{eq.subh}
\[
     \phi(\la_0) \,\le\, \frac{1}{2\pi}\int_0^{2\pi}
                         \phi(\la_0 + r e^{i\theta}) d \theta .
\]
%\en

A well known result from operator theory establishes
that for any bounded linear operator $\,T$ on a Hilbert space,
\[
     \|T\| \,=\, \sup_{\phi,\psi\ne 0} \frac{\textup{Re}
     \langle T \phi, \psi \rangle }{\|\phi\|\, \|\psi\|} \, .
\]
If the Hilbert space is finite dimensional, it is easy to see that
the supremum is attained. By virtue of (\ref{eq.l1}), $\,s_n(\la)
= \| P(\la)^{-1} \|^{-1}$. Thus,
\eq \label{eq.1}
    s_n(\la) \,=\, \left [ \max_{u,v\ne 0}
                   \frac{\textup{Re}\langle P(\la)^{-1}u,\,v\rangle }
                   {\|u\|\,\|v\|}  \right ]^{-1} .
\en Now note that $s_n(\la)^{-1}$ is continuous on $U$, and let
$\,\la_0 \in U\,$ and $\,u_0,v_0\in \CC^n\,$ be such that
\[
  s_n(\la_0)^{-1} \,=\, \frac{\textup{Re}\langle P(\la_0)^{-1}u_0,\,v_0\rangle }
                        {\|u_0\|\,\|v_0\|} \, .
\]
The function $\langle P(\la)^{-1}u_0,v_0\rangle $ is analytic on
$U$ and so the real function
\eq \label{eq.2}
  h(\la) \,:=\, \frac{\textup{Re}\langle P(\la)^{-1}u_0,\,v_0\rangle }
                {\|u_0\|\,\|v_0\|}
\en is harmonic on $U$. Furthermore, it follows from (\ref{eq.1})
and (\ref{eq.2}) that $\,h(\la) \le s_n(\la)^{-1}\,$ on $U$.
Consequently,
\[
    s_n(\la_0)^{-1} \,=\,  h(\la_0)
                    \,=\,  \frac{1}{2\pi} \int_0^{2\pi}
                           h(\la_0+re^{i\theta}) d \theta
                    \,\le\,\frac{1}{2\pi} \int_0^{2\pi}
                           s_n(\la_0+re^{i\theta})^{-1}d \theta ,
\]
and the result follows. \qed

\medskip

An important characteristic of subharmonic functions is the fact
that they satisfy the maximum principle. Therefore, the only local
minima of $s_n(\lambda)$ are those $\la\in\sigma(P)$.

The subharmonicity of $s_n(\la)^{-1}$ has been considered
recently by various authors. In \cite{BD}, Boyd and Desoer
discuss this property in the  context of linear control systems.
Concrete applications of this theorem may be found
in \cite{G} for the linear case, and in \cite{B,D2} for the
quadratic case. In \cite{B}, the result is applied in support of a
certain novel procedure for finding eigenvalues of self-adjoint
operators in infinite dimensional Hilbert spaces.

\begin{cor} \label{t2}
For all $\,\eps>0$, every connected component of $\lep$ has
non-empty interior.
\end{cor}
\noindent {\bf Proof.} Suppose, on the contrary, that $\GG$ is a
connected component of $\lep$ with empty interior. Since
$w(|\la|)/s_n(\la)$ is subharmonic, $\min_{\la\in\GG}
[s_n(\la)/w(|\la|)]$ should be attained at all points of $\GG$.
Thus, necessarily, $\GG$ should be a single point and in fact one
of the eigenvalues of $P(\la)$. The continuity of $s_n(\la)$ and
the fact that $\,w(0)=w_0>0\,$ ensure that this is not possible.
\qed

\medskip

In general, $s_n(\la)$ itself is not a subharmonic function as it
does not satisfy the maximum principle (a concrete example may be
found at the end of this section). However, as we will see next,
$s_n(x)$ is locally regular.

First consider the nonnegative
eigenvalue functions generated on $\CC$ by the matrix function
$P(\la)^*P(\la)$, say $\,S_1(\la),S_2(\la),\ldots,S_n(\la)$. They
can be organised in such a way that they have a strong smoothness
property.

\begin{lem} \label{lem.an}
For any given analytic curve $\,\zeta : \RR\rightarrow\CC,\,$ the
eigenvalues of $P(\la)^*P(\la)$ can be arranged in such way that,
for all $j$, $S_j(\zeta(t))$ are real analytic functions of
$\,t\in\RR$.

Furthermore, if $\,s_n(\la)=\min_j(S_j(\la))^{1/2}\,$ is a
non-zero simple singular value of $P(\la)$ and
$\,u_{\la},v_{\la}\,$ are associated left and right singular
vectors, respectively, then (writing $\,\la =x+iy$) $s_n(\cdot)$
is a real analytic function in a neighbourhood of $\lambda$ and
\eq \label{eq.grad}
       \nabla s_n(x+iy) \,=\,
       \left ( \textup{Re} \left ( u_{\la}^* \,
       \frac{\partial P(x+iy)}{\partial x }
       \, v_{\la} \right ) , \,
       \textup{Re} \left ( u_{\la}^* \,
       \frac{\partial P(x+iy)}{\partial y}
       \, v_{\la} \right )  \right ) .
\en
\end{lem}

The first statement follows from  Theorem S6.3 of \cite{GLR} (see
also Theorem II-6.1 of \cite{Ka}). For the second and third, see
\cite{S}, for example.

We can interpret the first part of this lemma pictorially in the
following manner. For $\,t\in\RR ,\,$ the graphs of
$S_j(\zeta(t))$ ($j=1,2,\dots ,n$) are smooth and might cross
each other. At a crossing point, the graph of the corresponding
singular value $s_k(\zeta(t))$ is continuous but it changes from
one smooth curve  to another with a possible jump in the
derivative (see \S II-6.4 of \cite{Ka}).

We may also consider regularity properties of $s_n(\lambda)$
as a function defined on the complex plane.
In this case, some rudimentary ideas from algebraic
geometry assist in
discussing the $n$ surfaces in $\RR^3$ which are (in general)
generated by the singular values. (Where possible, the terminology
of Kendig \cite{K} is followed).
Write $\,\la\in\CC\,$ in real and imaginary parts; $\,\la =
x + i y ,\,$ and define $n$ subsets of $\RR^3$:
\[
    \Sigma_j \,:=\, \left ( x , y , S_j(x+iy) \right )
             \; ; \;\; j = 1 , 2 , \ldots , n.
\]

\begin{propos} \label{prop.av}
The union $\;\bigcup_{j=1}^n\Sigma_j\,$ is a real algebraic
variety.
\end{propos}

\noindent {\bf Proof.} Define the function
\eq\label{eq.po}
   d(x,y,S)\,:=\,\det\left( I\,S - P(x+iy)^*P(x+iy)\right)
   \; ; \;\;  x , y , S \in \RR .
\en
Since the matrix $P(x+iy)^* P(x+iy)$ is hermitian, $d(x,y,S)$
is a polynomial in $\,x,y,S\,$ with real coefficients, and since
\[
   \bigcup_{j=1}^n\Sigma_j \,=\, \left \{ (x,y,S) \in \RR^3 :\,
                                      d(x,y,S) = 0 \right \} ,
\]
the result follows. \qed

\medskip

In spite of this proposition and Lemma~\ref{lem.an}, the
existence of an arrangement of the eigenvalues of $P(\la)^\ast
P(\la)$ such that the $n$ surfaces $\,\Sigma_j \in\RR^3\,$ are
smooth {\em everywhere} is not guaranteed in general. Consider
the following example. For the linear matrix polynomial
$\,P(\la)=I\la-A,\,$ where
\[
   A \,=\, \mat{ccc} 3/4 & 1 & 1 \\
                     0 & 5/4 & 1 \\
                     0 & 0 & -3/4   \rix ,
\]
$\Sigma_1\cup\Sigma_2\cup\Sigma_3$ has a conic double point at
$(0,0,5/16)$. Therefore, no arrangement of the singular values
exists ensuring $\,\Sigma_{1}$, $\Sigma_{2}\,$ and
$\,\Sigma_{3}\,$ are simultaneously smooth at $\,\la =0$. Moreover, in
this example, $\,s_3(0)=s_2(0)=\sqrt{5/16},\,$ so note that the
hypothesis of non-degeneracy of the fundamental singular value in
the second part of Lemma~\ref{lem.an} is essential.

For linear polynomials, the occurrence of isolated singularities
in $\bigcup_{j=1}^n\Sigma_j$ is rare. In the above example the
matrix $A$ had to be carefully crafted to allow the conic double
point around the origin. Any slight change in the coefficients of
$A$ would eliminate this degeneracy.

The following useful proposition is an immediate consequence of
Lemma~\ref{lem.an}.

\begin{propos} \label{prop.eq} If $\,S_j(\la)=S_k(\la)\,$
for $\,j\ne k\,$ and for all $\la$ in a non-empty open set
$\mathcal{O}$, then $\,\mathcal{O}=\CC$.
\end{propos}

Thus, different surfaces $\Sigma_j$ can intersect only in sets of
topological dimension at most one.

%%%%%%%%%%%%%%%%%%%%%%%%%%%%%%%%%%%%%%%%%%%%%%%%%%%%%%%%%%%%
\section{The pseudospectrum and its boundary} \label{sec.4}
%%%%%%%%%%%%%%%%%%%%%%%%%%%%%%%%%%%%%%%%%%%%%%%%%%%%%%%%%%%%

Now we put these ideas into the context of the study of
pseudospectra.

A fundamentally important case is that in which $\eps >0$ is so
small that $\lep$ consists of ``small'' disconnected components,
each one containing a single (possibly multiple) eigenvalue of
$P(\la)$. As $\eps$ is increased from zero, these components
enlarge, collide and eventually intersect in various ways so that
the boundary of $\lep$, say $\partial\lep$, becomes more complex.
In an earlier paper \cite{LP} two of the present authors studied
some basic properties of $\lep$ and $\partial\lep$ in support of a
curve-tracing algorithm for plotting $\partial\lep$.

Let
\begin{equation} \label{e20}
  F_{\eps}(x,y) \, \equiv \, F_{\eps}(x + i y)
    \, := \, s_n ( x + i  y ) - \eps \, w( | x + i y | )
                \; ; \;\;  x , y \in \RR.
\end{equation}
Since this function is continuous in $\,\la = x+iy \in\CC,\,$ it
follows from (\ref{eq.th}) that
\eq \label{eq.bo}
   \partial\lep \,\subseteq\,\{\la\in\CC :\,F_\eps(\la)=0 \}.
\en
Moreover, as long as $s_n(\la)^2$ is a {\em simple
non-vanishing} eigenvalue of $P(\la)^*P(\la)$, differentiation in
the direction of the boundary will be well-defined as a
consequence of Lemma \ref{lem.an}. However, when $\,s_{n-1}(\la) =
s_n(\la),\,$ this smoothness of the boundary may be lost.
Hence our interest in the set of $\,\la\in\CC\,$ for which $s_n(\la)$
is multiple; curve tracing algorithms are prone to fail around
these points, as the directional derivatives along
$\partial\lep$ may not be well-defined.

Even though it is quite rare\footnote{This is a rather delicate
point, and it seems to have been missed in the work of several
preceding authors as in \cite{LP} and \cite{TH}. In particular,
Corollary 4.3 of \cite{AB} seems to be false as it stands. On the
other hand, this fact seems to have little, if any impact on the
design of algorithms.}, in general, the right side of
(\ref{eq.bo}) might include points in the interior of $\lep$.
This can be observed as a consequence of either of the two
unlikely events:
\begin{description}
\item[(i)] the surface $s_n(\la)$ having a local
           (but not global) maximum,
\item[(ii)] at least three multiple sheets of $\bigcup_{j=1}^n
            \Sigma_j$ intersecting in a single point.
\end{description}

Demmel's matrix
\[
    A=\left(\begin{array}{ccc} -1 & -b & -b^2 \\ 0 & -1 & -b \\
0 & 0 & -1     \end{array}\right)
\] 
with $b>\!>1$ illustrates (i) for the standard eigenvalue problem with 
$w(x)=1$.
Indeed if $b=100$ and $P(\lambda)=(\lambda I-A)$, $s_3(\la)$
has a local maximum at $\la=0$, cf. \cite{demmel}.  

Higher order examples typifying (i) can also be easily constructed. 
Consider, for instance, the polynomial
$\,P(\la)=(\la^2-1)(\la^2-i)\,$ in $\CC$ and the weight function
$\,w(x) = 4x^2 + 1$. The point $\,\la=0\,$ is a local maximum of
the function
\[
      \frac{s_1(\la)}{w(|\la|)} \,=\,
      \frac{|\la^2-1|\,|\la^2-i|}{4|\la|^2+1} \, ,
\]
which is smooth in $\,\CC \, \backslash \, \{\pm 1,\,\pm
i^{1/2}\}$. This may be verified by directly computing the
gradient and Hessian of this expression at $\,\la=0$. Thus, when
$\eps =1$ and $\la$ lies in a sufficiently small neighbourhood
$\N$ of the origin we have $\,s_1(\la) \le w(|\la|),\,$ so that
$\,\N\in \Lambda_1(P)\,$ and $\,0\notin \partial\Lambda_1(P)$.
However, for $\eps =1$,
\[
    F_{\eps}(0) \,=\, F_1(0) \,=\, s_1(0)-w(0) \,=\, 0 ,
\]
so in this case the inclusion of (\ref{eq.bo}) is proper.

To confirm (ii), recall Example 3.5 of \cite{AB}: for $\,w(x)=1\,$
and any $\,\eps >0,\,$ the $\eps$-pseudospectrum of $\,P(\la) =
\textup{diag}\{\la -1,\la+1,\la-i,\la+i\}\,$ is the union of four
closed discs with centres at the eigenvalues $\,1,\,-1,\,i,\,-i\,$
and radii equal to $\eps$. Thus, for $\,\eps=1$, the origin lies
in the set $\,\{ \la \in \CC :\,F_1 (\la)=0 \}\,$ but it is an
interior point of $\Lambda_1(P)$.

The next result shows that $\partial\lep$ is made up of
algebraic curves. This is a comforting property in the sense that
the number of difficult points, such as cusps or
self-intersections, is limited. (See Proposition 6.2.10 of
\cite{BK} for an explicit statement of this kind.)

\begin{thm} \label{p7}
Let $\eps>0$ and assume that $\,\lep\ne\CC$. Then the boundary of
$\lep$ lies on an algebraic curve. In particular, $\partial\lep$
is a piecewise $C^\infty$ curve, it has at most a finite number
of singularities where the tangent fails to exist, and it
intersects itself only at a finite number of points.
\end{thm}
\noindent{\bf Proof.}
We first show that $\partial\lep$ lies on an algebraic curve.
Recall the function $d(x,y,S)$ defined by
(\ref{eq.po}) and observe that $\partial\lep$ lies on the level
set
\begin{eqnarray*}
  \LL_1 &=& \left \{ x+iy :\, x,y \in \RR, \; \eps \, w(|x+iy|)
            \;\,\mbox{is a singular value of}\,\; P(x+iy) \right \} \\
        &=& \left \{ x+iy :\, x,y \in \RR, \;
            d(x,y,\eps^2 w(|x+iy|)^2) = 0 \right \} .
\end{eqnarray*}
The function $d(x,y,\eps^2 w(|x+iy|)^2)$ can be written in the form
\[
  d(x,y,\eps^2 w(|x+iy|)^2) \,=\, \sqrt{x^2+y^2}\, p(x,y) + q(x,y),
\]
where $p(x,y)$ and $q(x,y)$ are real polynomials in $\,x,y\in\RR$.
Thus,
\[
   \LL_1 \,=\, \left \{ x+iy : \, x,y \in \RR, \;
               \sqrt{x^2+y^2}\, p(x,y) + q(x,y) = 0 \right \} .
\]

If $w(x)$ is an even function, then $p(x,y)$ is identically zero
and either $\LL_1$ is an algebraic curve or it coincides with the
complex plane. Suppose $w(x)$ is not an even function. Then
$\LL_1$ is a subset of the level set
\[
   \LL_2 \,:=\, \left \{ x+iy : \, x,y \in \RR , \;
                 (x^2+y^2) p(x,y)^2 - q(x,y)^2 = 0 \right \} ,
\]
which is also an algebraic curve when it does not coincide with the
complex plane.

Next we show that $\,\LL_2=\CC\,$ only if $\,\LL_1=\CC$. Thus, if
$\LL_2=\CC$ and $p(x,y)$, $q(x,y)$ are not identically zero, then
\[
     (x^2+y^2) \, p(x,y)^2 \,=\, q(x,y)^2
              \;\;\, \mbox{for all}\;\; x , y \in \RR ,
\]
where the order of the (irreducible) factor $\,x^2 + y^2\,$ in the
left hand side is odd and the order of the same factor on the
right (if any) is even. This is a contradiction. Hence, if
$\,\LL_2=\CC,\,$ then $p(x,y)$ and $q(x,y)$ are identically zero,
and consequently, $\,\LL_1=\CC$.

Since $\,\LL_1 \subseteq \lep$ and by hypothesis $\,\lep\ne
\CC,\,$ both $\,\LL_1 ,\LL_2\ne\CC\,$ and so $\LL_2$ is an
algebraic curve.  This completes the first part of the theorem.

For the second part, note that, as $s_n(\la)$ is continuous in
$\la\in\mathbb{C}$, $\partial\lep$ is a union of continuous
curves. From the above considerations it follows that $\LL_1$ is
a piecewise $C^\infty$ curve and it has finitely many
singularities. Then, since
\[
  \partial \lep \,\subseteq\, \LL_1 \,\subseteq\, \lep ,
\]
we can actually decompose $\,\LL_1 = \bigcup_{k=1}^h \gamma_k
,\,$ where $\gamma_k$ ($k=1,2,\dots , h$) are suitable smooth
curves with the following property: $\, \gamma_k \subseteq
\partial \lep\,$ for all $\,1 \leq k \leq j ,\,$ and
$\,\gamma_k \subseteq \lep \setminus \partial \lep\,$ for all
$\,j < k \leq h$. Thus, $\,\partial\lep = \bigcup_{k=1}^j
\gamma_k\,$ as needed. \qed

\medskip

Note that for the standard eigenvalue problem, $\,w(x)=1\,$ is
an even function. In this case, the above result appears in
the work of Alam and Bora \cite{AB}.

The following technical statements will be useful subsequently.
The first one follows immediately from (\ref{eq.bo}).

\begin{lem} \label{lem.bs}
If $\,0\leq\delta <\eps,\,$ then $\,\partial\Lambda_{\delta}(P)
\subset\lep\,$ and $\,\Lambda_{\delta}(P)\cap\partial\lep =
\emptyset$.
\end{lem}

In particular, note that $\,\sigma(P)\cap\partial\lep =
\emptyset\,$ for any $\,\eps>0$.

With $P(\la)$ as in (\ref{eq.poly}), consider a perturbed matrix
polynomial $\,Q(\la)$ of the form (\ref{e1}). It follows
from the definition (\ref{eq.b}) that $\,Q(\la) \in \partial
\B(P,\eps,w)\,$ if and only if $\,\|\Delta_j\| \le \eps \, w_j\,$
for each $j$ and equality holds for at least one $j$.
Now consider matrix polynomials in the interior of $\B(P,\eps,w)$;
$\interior[\B(P,\eps,w)]$. It is easily seen that $\,Q(\la) \in
\interior[\B(P,\eps,w)]\,$ if and only if
\[
     \|\Delta_j\| \,<\, \eps \, w_j \;\;
     \textup{whenever} \;\, w_j > 0 , \;\, \mbox{and}
\]
\[
     \Delta_j \,=\, 0  \;\; \textup{whenever} \;\, w_j = 0 .
\]

\begin{lem} \label{lem.con} If $\,\mu\in\partial\lep ,\,$ then
for any perturbation $\,Q(\la)\in\B(P,\eps,w)\,$ such that
$\,\mu\in\sigma(Q),\,$ $\,Q(\la)\in\partial\B(P,\eps,w)$.
\end{lem}

\noindent {\bf Proof.} Let $\,\mu\in\partial\lep$. It suffices to
show that if $\mu\in\sigma(Q)$ for a $\,Q(\la)\in
\mathcal{B}(P,\eps,w),\,$ then $\,\|\Delta_j\| = \eps\, w_j\,$ for
some $\,j=0,1,\ldots ,m$. Indeed, if we assume the converse
statement, $\,\|\Delta_j\| < \eps\, w_j \,$ for all $j$, then
$\,Q(\la)\in\mathcal{B}(P,\tilde\eps,w)$ for some $\tilde\eps <
\eps$. But since $\,\mu\in\sigma(Q),\,$ we have $\,\mu\in
\Lambda_{\tilde\eps}(P),\,$ which contradicts Lemma \ref{lem.bs}.
Thus, the desired assertion holds. \qed

%%%%%%%%%%%%%%%%%%%%%%%%%%%%%%%%%%%%%%%%%%%%%%%%%%%%%%%%%%%%
\section{The fault lines} \label{s4}
%%%%%%%%%%%%%%%%%%%%%%%%%%%%%%%%%%%%%%%%%%%%%%%%%%%%%%%%%%%%

Differentiability along $\partial \lep$, the boundary of the
pseudospectrum, is possible as long as the gradient of
$\,s_n(\la)-\eps\, w(|\la|)\,$ exists and does not vanish. The
only place where $w(|\la|)$ might fail to have a derivative is the
origin. If the minimal singular value, $s_n(\la)$, has
multiplicity one, then $s_n(\la)$ is smooth in a neighbourhood of
$\la$. Thus, the study of those points where differentiability is
lost, apart from $\,\la = 0,\,$ is confined to the region of the plane
where the sheet of $\bigcup_{j=1}^n\Sigma_j$ corresponding to $s_n(\la)$,
meets the one corresponding to  $s_{n-1}(\la)$. This motivates the
following definition.

The rather involved indexing introduced next is required 
to handle multiple eigenvalues.
Below we always assume that the eigenvalues of $P(\la)^\ast
P(\la)$ are ordered so that $S_j(\la)=s_j(\la)^2$ for all
$j=1,2,\dots,n$. Let \linebreak $p:\{1,2,\dots,n\}\longrightarrow
\{1,2,\dots,n\}$ (usually not onto) satisfying the following
properties:
\begin{description}
\item[(a)] $\Sigma_j=\Sigma_{p(j)}$,
\item[(b)]
$\bigcup_{j=1}^n \Sigma_j=\bigcup_{j=1}^n \Sigma_{p(j)},\,$ and
\item[(c)] $\Sigma_{p(j)}=\Sigma_{p(k)}\,$ if and only if $\,p(j)=p(k)$.
\end{description}
The map $p$ is a choice of the indices of those, and only those,
different $\Sigma_j$. Let $\,c_1:=\max\{p(j)\}_{j=1}^n\,$ and
$\,c_2:=\max[\{p(j)\}_{j=1}^n \setminus \{ c_1 \}]$. We define the
set
\[
  \mathcal{F}_P \,:=\, \{\la \in \mathbb{C}\,:\,s_{c_1}(\la)
                          = s_{c_2}(\la) \} .
\]

By virtue of Proposition \ref{prop.eq}, $\F_P$ has empty interior.
Furthermore, if all eigenvalues of $P(\la)$ have 
geometric multiplicity equal to $1$, then
\[
  \mathcal{F}_P \,=\, \{\la\in \CC:s_n (\la)=s_{n-1}(\la)\}.
\]

\begin{propos}
If all the eigenvalues of $P(\la)$ have geometric multiplicity
equal to $1$, then either $\,\F_P=\emptyset\,$ or $\,\F_P$ lies
on an algebraic curve (including the possibility of isolated
points).
\end{propos}

\noindent {\bf Proof.}  Let
\[
  \hat{\F} \,=\, \left\{\la\in\CC :\, s_j(\la) = s_k(\la),\;
                     j \ne k \right \}
\]
so that $\,\F_P\subset\hat{\F}$. This set is the locus of all
points $\,(x,y)\in \RR^2\,$ such that the discriminant of the real
polynomial in $S$ defined by (\ref{eq.po}) is zero.
The hypothesis ensures that $\,\hat{\F}\ne\CC,\,$ and thus, either
$\,\hat{\F}=\emptyset\,$ or $\hat{\F}$ is an algebraic curve. The
result follows just because $\F_P$ is a subset of $\hat{\F}$. \qed

\medskip

In particular, $\mathcal{F}_P$ might include straight lines, 
single points, the empty set, or be a complicated set such as a
Voronoi diagram (see Example 1 below).

Borrowing a geological term, we call the set $\F_P$ the set of
{\em fault points} of $P(\la)$. In general, $\F_P$ will be made up of
{\em fault lines}. The explicit determination of the fault lines
of $P(\la)$ requires computations with determinants and
discriminants, and is therefore unrealistic. However, the
following considerations demonstrate the role that the fault
lines frequently play in the study of pseudospectra.

Let $F_\eps(x,y)$ be as in (\ref{e20}). As mentioned above, apart
from $\,\la=0,\,$ if $\nabla F(x,y)$ does not exist, then $\,x+iy
\in\F_P$. At these points, the curve $\partial\lep$ will typically
fail to have a tangent line. There are other points where the
tangent line will be undefined, those where $\,\nabla F(x,y)=0$.
In this case, there is a saddle point in the minimal singular
value surface. These may or may not lie on $\mathcal{F}_P$
(see Section \ref{section:multiple}).

Example 3 below illustrates a case in which $\mathcal{F}_P$ is a
singleton. In Example~4, $\mathcal{F}_P$ is empty but there is,
nevertheless, a point at which $\partial\lep$ has no tangent. In
Examples 1, 2, and 5, $\mathcal{F}_P$ is, indeed, made up of
fault {\em lines}.

\medskip

\noindent \textbf{Example 1} Let $A$ be an $\,n \times n\,$ normal
matrix with eigenvalues $\,\{\la_j\}_{j=1}^n$. Then the
fault lines of $\,P(\la)=I \la -A\,$ (i.e., for the standard
eigenvalue problem) form the Voronoi diagram defined by
$\,\{\la_j\}_{j=1}^n\,$ (i.e., the boundary of their Dirichlet
tessellation). \qed

\medskip

\noindent \textbf{Example 2} Naive experiments with diagonal
matrix polynomials provide an insight on the possible structure of
individual fault lines. For instance, let $\,P(\la) =
\textup{diag} \{ \la^2-2\la ,(a-\la)(\la+2)\}\,$ and set
$\,w(x)=1$.
%%%%%%%%%%%%%%%%%%%%%%%%%%%%%%%%%%%%%
\begin{figure}
\epsfig{file=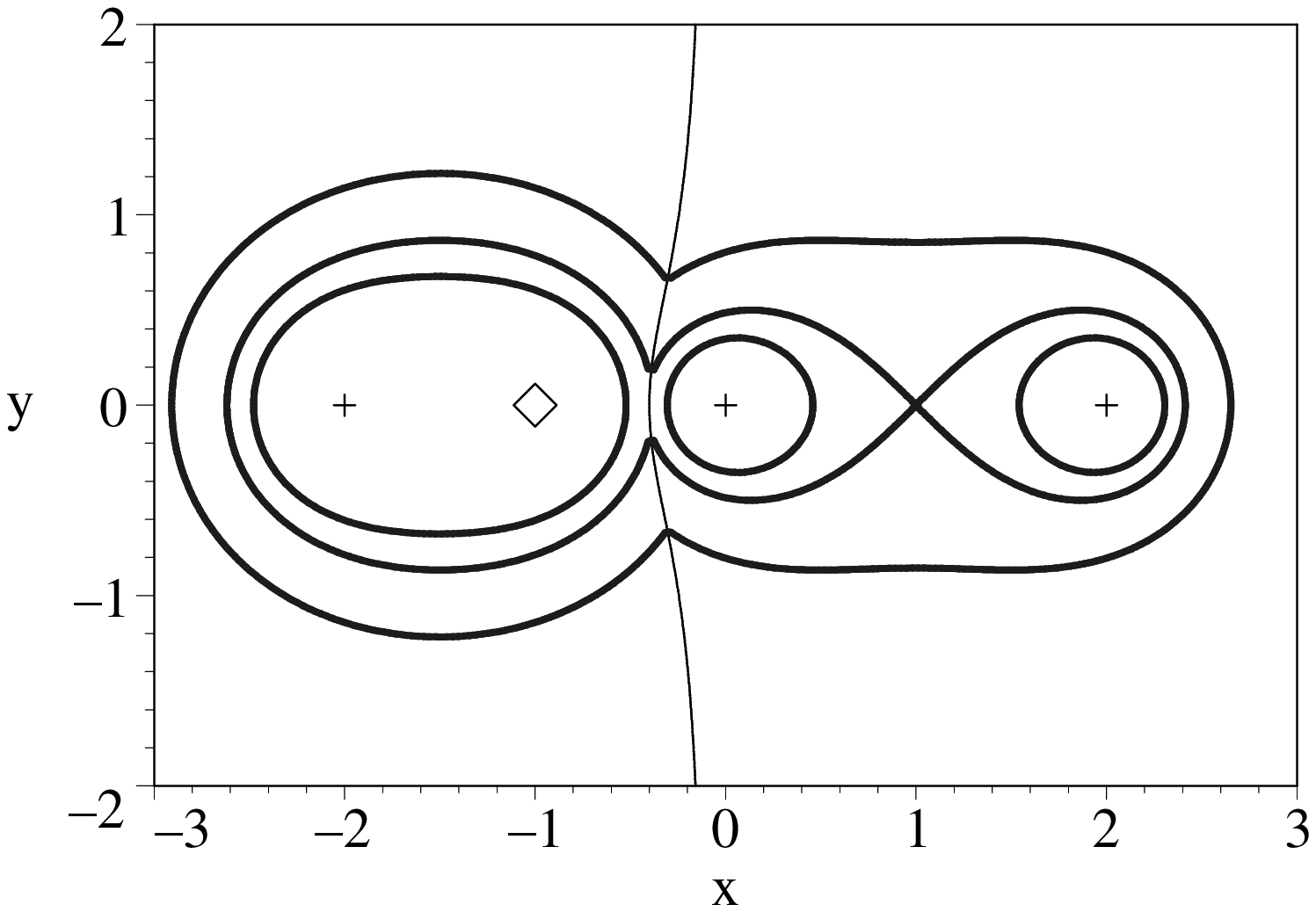, width=2.7in} %%%
\epsfig{file=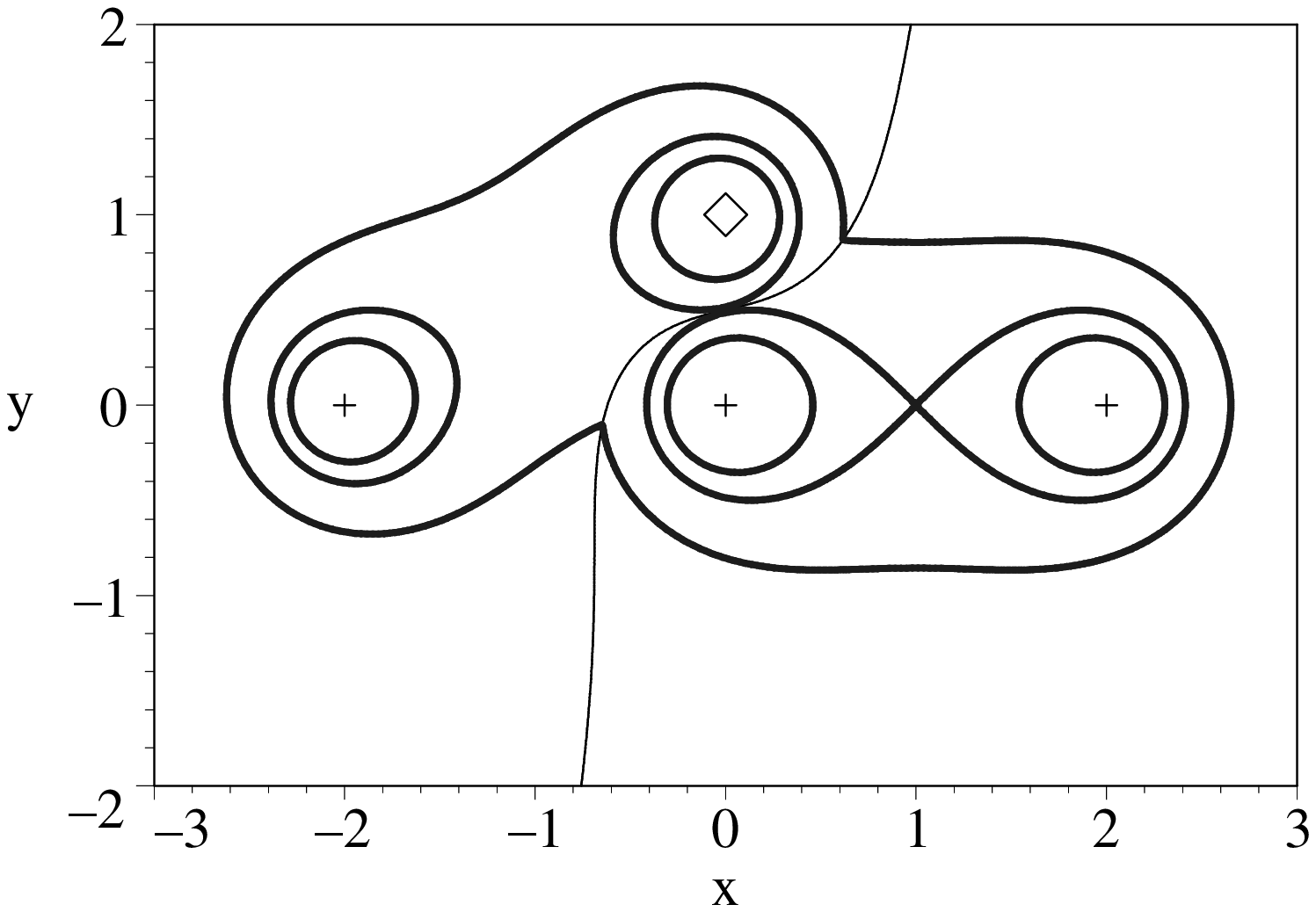, width=2.7in} %%%
\epsfig{file=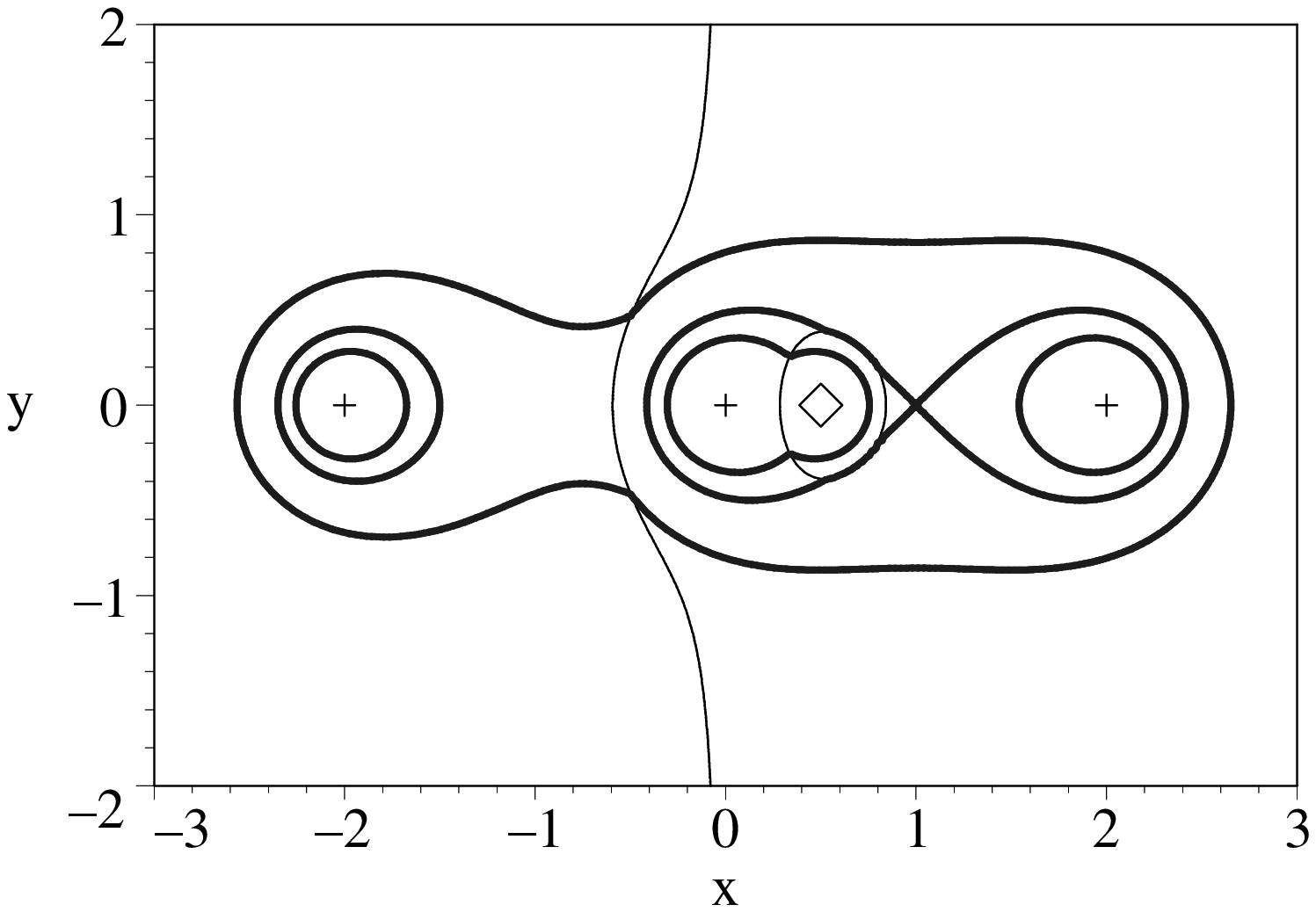, width=2.7in} %%%
\epsfig{file=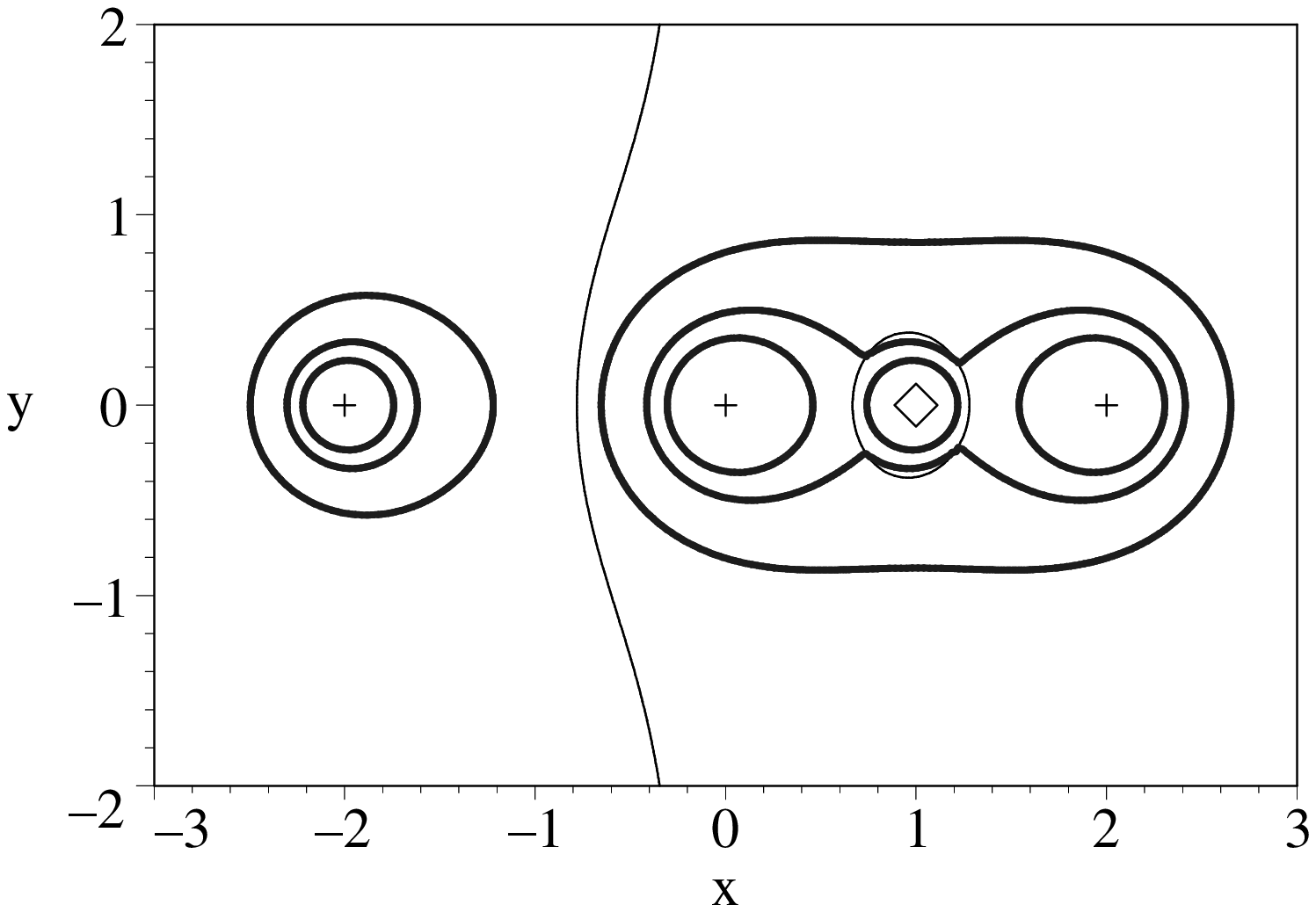, width=2.7in} %%%
\caption{\small The thin solid lines are $\mathcal{F}_P$. The thick
solid lines are $\partial\lep$. \label{f3}}
\end{figure}
%%%%%%%%%%%%%%%%%%%%%%%%%%%%%%%%%%
In Figure \ref{f3}, we depict the evolution of the set
$\mathcal{F}_P$ and $\partial\lep$ ($\eps = 1/\sqrt{2},\, 1,\,
\sqrt{3}$) for $\,a = -1,\, i,\, 1/2,\, 1$. The fixed eigenvalues
of $P(\la)$ are plotted as ``+'' and the perturbed eigenvalue $a$
is marked with a ``$\diamond$''.

In general, an unbounded $\mathcal{F}_P$ appears to be more
likely to occur. Nonetheless this set can also contain a compact
smooth curve. In both of the lower figures, $\mathcal{F}_P$
consists of an unbounded curve, which is asymptotic to a vertical
line, and a closed compact curve on the right half plane around
the perturbed eigenvalue $a$. The curve $\partial\Lambda_1(P)$
has a self intersection at $\,\la=1\,$ for $\,a = -1,\, i,\,
1/2$. This can be shown from the fact that this part of the
pseudospectrum depends only on the first diagonal entry of
$P(\la)$. The self intersection disappears as soon as $a$ moves
sufficiently close to $1$. There are critical values of $a$,
where $\,1 \in \mathcal{F}_P$. Two of these critical values are
$\,a=2/3\,$ and $\,a=4/3$. \qed

\medskip

\begin{figure}
\hspace{.1in} \epsfig{file=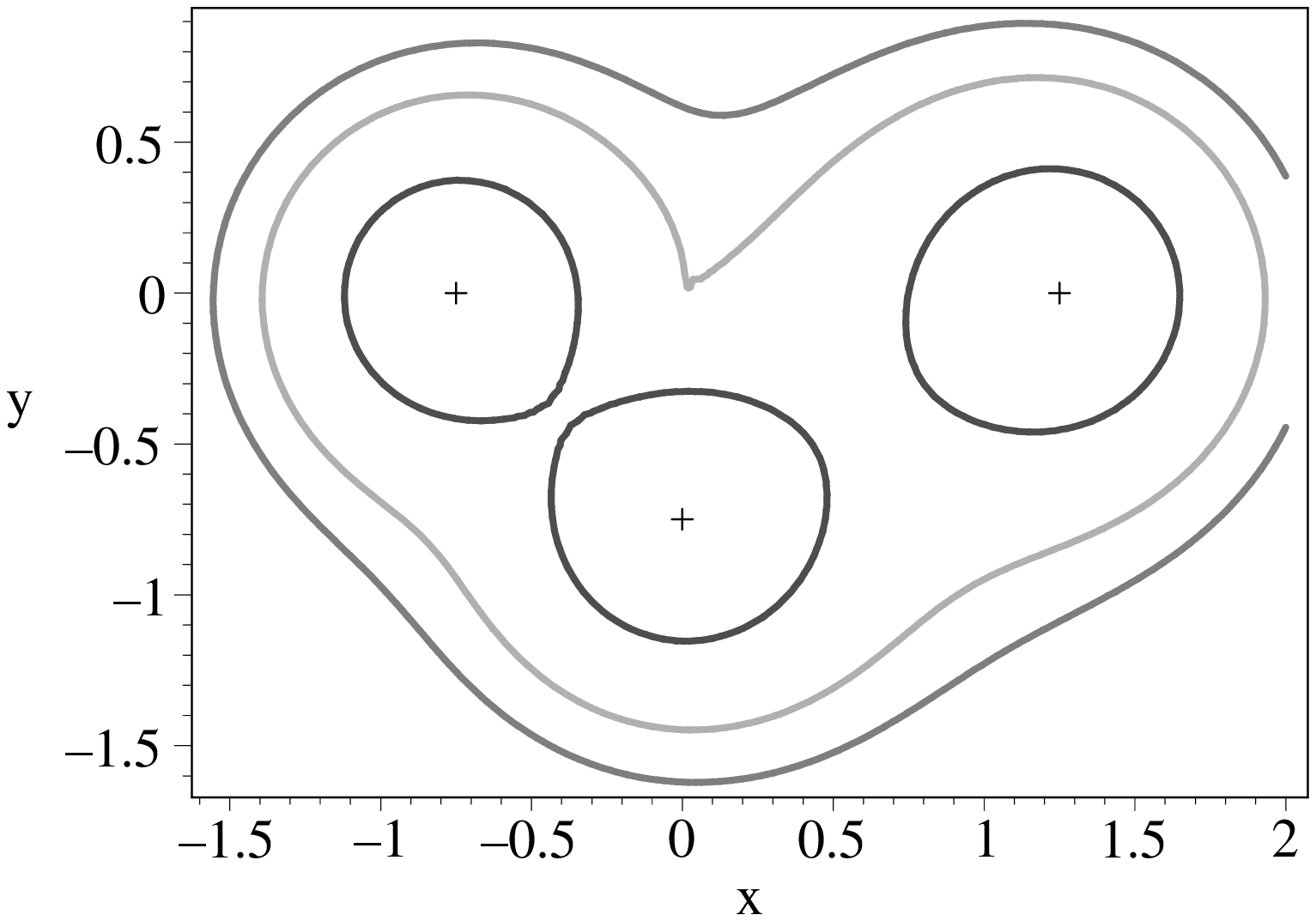, width=2.7in} %%%
\hspace{.4in}\epsfig{file=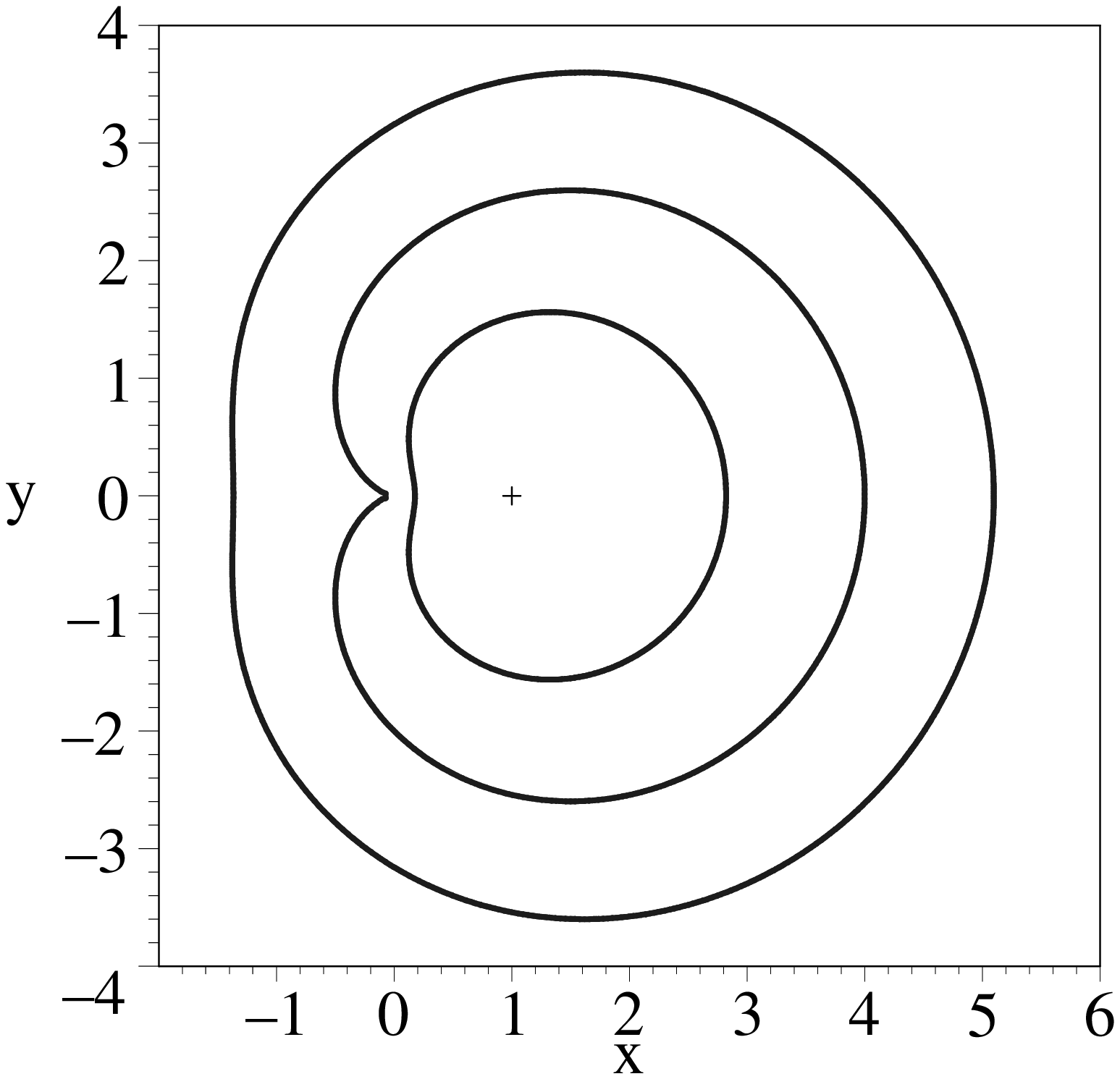, width=2in}    %%%
\caption{\small $\partial\lep$ has no tangent line at the origin
for just one $\eps$. \label{f4}}
\end{figure}

\noindent \textbf{Example 3} $\mathcal{F}_P$ can also be a singleton.
In the left part of Figure \ref{f4}, we depict
$\partial\lep$ for the linear matrix polynomial
\[
   P(\la) = \mat{ccc} \la + 3 i/4 & 1 & 1  \\
                           0 & \la-5/4 & 1 \\
                           0 & 0 & \la + 3/4 \rix ,
\]
the weight function $\,w(x)=1\,$ and $\,\eps^2 = 1/10,\, 5/16,\,
1/2$. The very special structure of this matrix polynomial ensures
that $\,\mathcal{F}_P=\{0\}$. The boundary of the pseudospectrum
does not have a tangent line at $\,\la=0\,$ when
$\,\eps=\sqrt{5/16}$. Compare with Example 4 below. \qed
%%%%%%%%%%%%%%%%%%%%%%%%%%%%%%%%%%%
%%%%%%%%%%%%%%%%%%%%%%%%%%%%%%%%%%%

\medskip

By construction, $\F_P$ is independent of  $w(x)$. Therefore, the
singularities occurring on $\partial \lep$ in places where the
gradient of (\ref{e20}) fails to exist, are, with the possible
exception of $\la=0$, independent of the chosen weights. In order
to illustrate this remarkable fact, we consider two more
examples.

\medskip

\noindent \textbf{Example 4} The set $\F_P$ might be empty but the
smoothness of $\partial\lep$ might be broken at $\,\la=0\,$ due to
the weight function. Indeed, let $\,n=1$, $\,P(\la)=(\la-1)^2\,$
and $\,w(x)=2x+1$. Then $\,s_1(x+iy)=(x-1)^2+y^2\,$ and
$\,\mathcal{F}_P = \emptyset$.

When $\,\eps=1$, $\,F_1(x,y) = x^2
+ y^2 - 2x - 2 \sqrt{x^2+y^2}$. Hence, $\,F_1(x,y)=0\,$ if and
only if
\[
   (x-1)^2+y^2 \,\ge\, 1  \;\; \mbox{and} \;\;
   y^4 + 2(x^2-2x-2) y^2 + (x^4-4x^3) \,=\, 0 .
\]
Thus, the curve $\partial \Lambda_1(P)$ has a parameterisation of
the form
\[
    y_\pm(x) \,=\, \pm \,\sqrt{2+2x-x^2-2\sqrt{2x+1}}
    \; ; \;\; -1/10 \leq x \leq 0
\]
in a neighbourhood of the origin. As $\,\partial_x y_+(0) < 0 \,$
and $\,\partial_x y_-(0) > 0,\,$ $0\in \partial \Lambda_1 (P)\,$
is a singularity of Lipschitz type. The boundaries of $\lep$ for
$\,\eps = 1/2,\, 1,\, 3/2 ,\,$ are drawn in the right part of
Figure \ref{f4}. \qed

\medskip

\begin{figure}
\epsfig{file=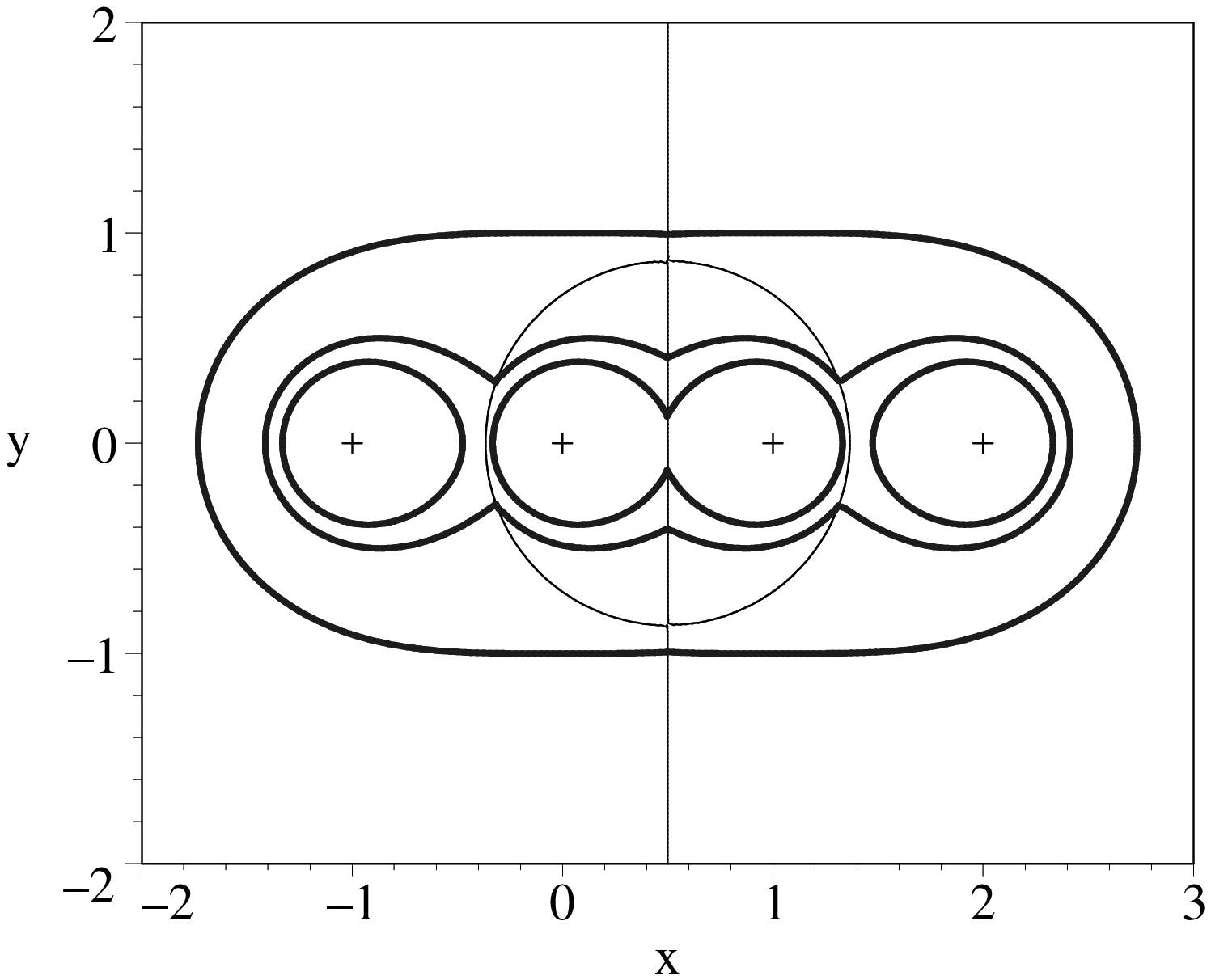, width=2.7in} %%%
\epsfig{file=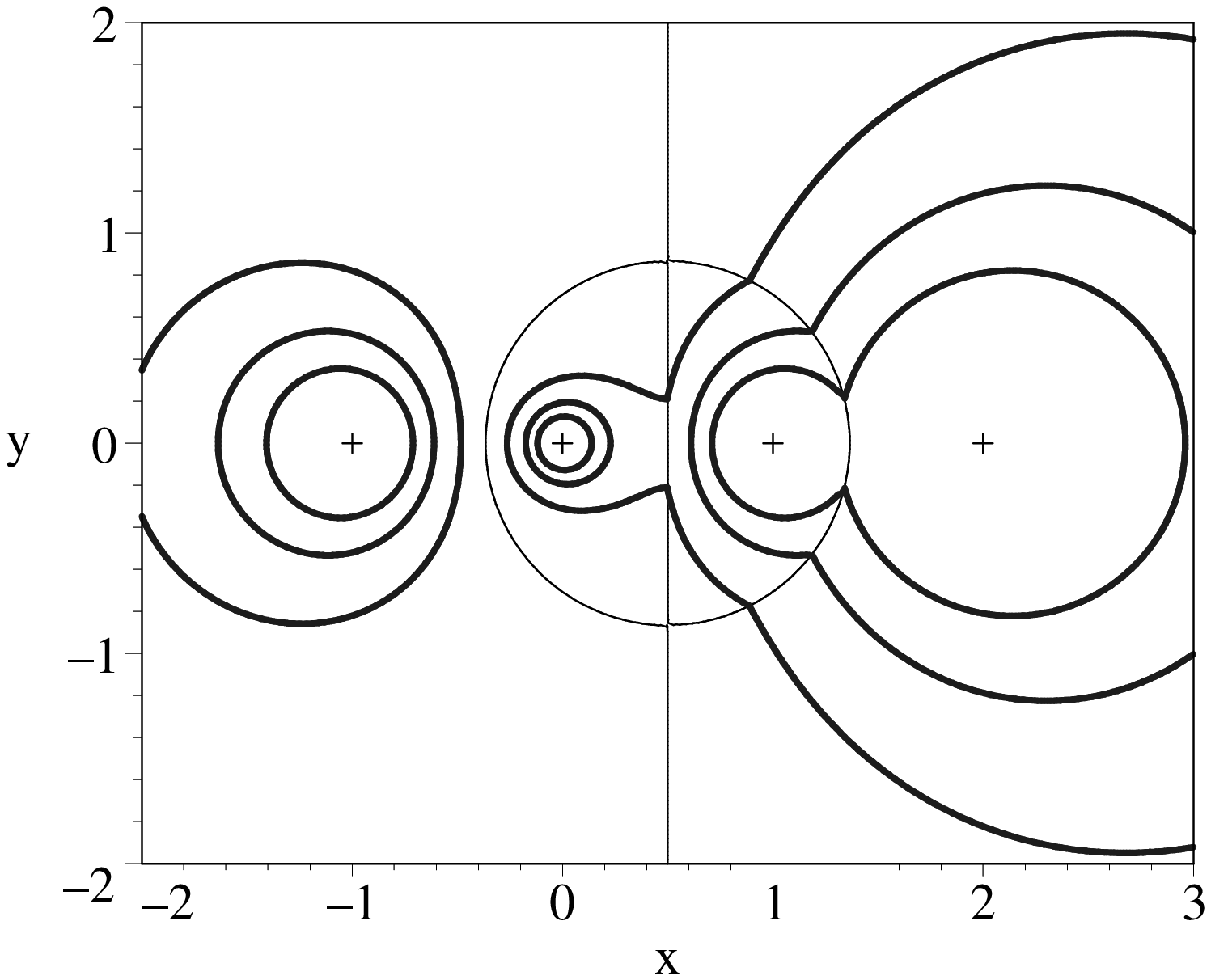, width=2.7in} %%%
\caption{\small The thin solid lines are $\mathcal{F}_P$. The thick
solid lines are $\partial\lep$. \label{f1}}
\end{figure}

\noindent \textbf{Example 5} Let $\,P(\la) = \textup{diag} \{
\la^2 - 1, \la^2 - 2 \la \}$. In Figure \ref{f1}, we depict
$\mathcal{F}_P$ and $\partial\lep$ for $\,w(x)=1\,$ and $\,\eps =
\sqrt{3/5},\, 1,\, 2\,$ (left), and for $\,w(x)=x^2+x+1\,$ and
$\,\eps^2 = 1/20,\, 1/10,\, 1/5\,$ (right). Here, $\mathcal{F}_P$
comprises a circle centred at $(1/2,0)$ and the line $\,x=1/2$.
As in the previous examples, ``+'' marks the locations of the
eigenvalues of $P(\la)$.   \qed

\medskip

All the above examples were designed in such a manner that both
the fault points and the boundaries of pseudospectra can be
constructed analytically either by hand or using algebraic
computer packages. We produced Figures~\ref{f3}, \ref{f4} and
\ref{f1} using commands provided in the standard distribution of
Maple.

%%%%%%%%%%%%%%%%%%%%%%%%%%%%%%%%%%%%%%%%%%%%%%%%%%%%%%%%%%%%%%%
\section{On the number of connected components} \label{number}
%%%%%%%%%%%%%%%%%%%%%%%%%%%%%%%%%%%%%%%%%%%%%%%%%%%%%%%%%%%%%%%

Consider an $\,n\times n\,$ matrix polynomial $P(\la)$ as in
(\ref{eq.poly}), a real $\,\eps>0,\,$ and a weight function
$w(\la)$ with $\,w(0)=w_0>0$. Theorem 2.3 of \cite{LP} will be
useful in the remainder of the paper. First we examine the case in which 
$\sigma(P)$ contains
multiple eigenvalues more carefully and without the
restriction of boundedness. A technical lemma will assist in the
argument.

\begin{lem}  \label{connection}
Suppose $A$ and $E$ are two $\,n \times n\,$ complex matrices such
that the determinants $\,\det A\,$ and $\,\det(A+E)\,$ are
nonzero. Then there is a continuous map $\,t \mapsto E(t) \in
\CC^{n \times n}$, $\,t \in [ 0 , 1 ] ,\,$ such that $\,E(0)=0$,
$\,E(1)=E,\,$ and
\[
     \det(A+E(t)) \,\ne\, 0 \;\; \mbox{and} \;\;
     \| E(t) \| \le  \| E \| \; ; \;\; t \in [ 0 , 1 ].
\]
\end{lem}

\noindent \textbf{Proof.} Since det$A \ne 0$ and det$(A+E)\ne 0$,
no eigenvalue of the pencil $A+tE$ can be equal to $0$ or $1$ (and
some may be infinite). So it may be assumed that det$(A+tE)$ has
$s$ real zeros in the interval $(0,1)$, where $\,0\le s \le n$.

If $\,s=0,\,$ then the continuous map $\,t \mapsto t E $, $\,t
\in [ 0 , 1 ] ,\,$ has the properties required by the lemma. If
$\,s \geq 1,\,$ then let $\,t_1 < t_2 < \cdots < t_s\,$ denote the
zeros of $\,\det(A+tE)\,$ in $\,(0,1)$. For any $\,t_j\,$ $(j=1,
2,\dots,s),\,$ the matrix $\,A+t_j E\,$ is singular and for
$\,\delta_j > 0\,$ sufficiently small, we have
\[
  \det [ A+(t_j+\ee^{i\theta}\delta_j) E ] \ne 0
     \;\; \mbox{ and } \;\;
  \| (t_j+\ee^{i\theta}\delta_j) E \| \le \| E \|
    \; ; \;\; \theta \in [ 0 , 2 \pi ] .
\]
In $\,[0,1],\,$ we replace each interval $\,[t_j-\delta_j,t_j+
\delta_j]\,$ with the circular arc
\[
 \mathcal{C}_j \,=\, \{ t_j - \ee^{i\theta} \delta_j :
                        \theta \in [ 0 , \pi ] \} ,
\]
and consider the continuous curve
\[
   \mathcal{S} = [ 0 , t_1 - \delta_1] \cup \mathcal{C}_1 \cup
   [ t_1 + \delta_1 , t_2 - \delta_2 ] \cup \mathcal{C}_2
   \cup \cdots \cup
   [ t_{s-1} + \delta_{s-1} , t_s - \delta_s ] \cup \mathcal{C}_s
   \cup [ t_s + \delta_s , 1 ]
\]
in the complex plane. For every continuous map $\,t \mapsto z(t)
\in \mathcal{S}$, $\,t\in [0,1],\,$ such that $\,z(0)=0$ and
$\,z(1)=1,\,$ the map $\,t\mapsto z(t) E \in \CC^{n \times n}$,
$\,t\in[0,1],\,$ has the required properties. \qed

%%%%%%%%%%%%%%%%%%%%%%%%%%%%%%%%%%%%%
%\medskip
%
%The main result of this section is:
%%%%%%%%%%%%%%%%%%%%%%%%%%%%%%%%%%%%%

\medskip

We are now ready to establish our main result on the number of
connected components of pseudospectra.
We should remark that, when $\lep$ is bounded, the following theorem
is a consequence of Theorem \ref{the.sa}.
Indeed, since $w(x)$ is a real polynomial,
$w(|\la|)$ is a subharmonic function in $\CC$ so, by Theorem
\ref{the.sa}, $s(\la)^{-1}w(|\la|)$ is subharmonic in
$\CC\setminus \sigma(P)$. If $\lep$ had a connected component 
where there is no
eigenvalue of $P(\la)$, then $s_n(\la)\,w(|\la|)^{-1}$ would have
a local minimum in this component, which is impossible according
to Theorem \ref{the.sa}.

\begin{thm} \label{thm.main3}
If the matrix polynomial $P(\la)$ has exactly $\,k$ $(\le n m)\,$
distinct eigenvalues (not necessarily simple), then for any
$\,\eps > 0,\,$ the pseudospectrum $\lep$ has at most $\,k$
connected components.
\end{thm}

\noindent \textbf{Proof.} If $\,\lep =\CC,\,$ then there is
nothing to prove. So assume that $\,\lep\ne\CC,\,$ and consider a
perturbation
\[
   Q (\la) \,=\, (P_m + \Delta_m) \la^m +
                 \cdots + (P_1 + \Delta_1) \la + P_0 + \Delta_0
\]
in $\mathcal{B}(P,\eps,w)$ with $\,\det (P_m + \Delta_m) \ne 0$.
By Lemma \ref{connection}, there is a continuous map $\,t \mapsto
\Delta_m (t) \in \CC^{n \times n}$, $t \in [ 0 , 1 ]$, such that
$\,\Delta_m(0)=0$, $\,\Delta_m(1)=\Delta_m,\,$ and
\[
 \det( P_m + \Delta_m(t) ) \,\ne\, 0 \;\; \mbox{and} \;\;
 \| \Delta_m(t) \| \le \| \Delta_m \| \; ; \;\; t \in [ 0 , 1 ].
\]
Hence, every member of the family
\[
   Q_t(\la) \,=\, ( P_m + \Delta_m(t) ) \la^m +
                  \cdots + (P_1 +  t \Delta_1) \la +
                  P_0 + t \Delta_0 \; ; \;\; t \in [0,1]
\]
has exactly $nm$ eigenvalues, counting multiplicities. Moreover,
all $Q_t(\la)$ $(t\in [0,1])$ belong to $\mathcal{B}(P,\eps,w)$.
Their eigenvalues lie in $\lep$ and trace continuous curves from
the eigenvalues of $P(\la)$ $( = Q_0(\la))\,$ to the eigenvalues
of $Q(\la)$ $( = Q_1(\la))$. Thus, as in the proof of Theorem 2.3
of \cite{LP}, the set
\[
   \Lambda_0 \,=\, \left \{ \mu \in \CC :\, \det Q(\mu) = 0 ,\;
   Q(\la)\in\mathcal{B}(P,\eps,w) ,\; \det (P_m + \Delta_m)
   \ne 0 \right \}
\]
has at most $k$ connected components determined by the $k$
distinct eigenvalues of $P(\la)$.

Now let $\la_0$ be an interior point of $\lep$, and let $\,R(\la)
= \sum_{j=0}^{m} R_j \la^j\,$ be a perturbation in
$\mathcal{B}(P,\eps,w)$ with $\,\det R_m = 0,\,$ such that
$\,\la_0\in \sigma(R)$. Since $\,\lep\ne\CC,\,$ $R(\la)$ has less
than $nm$ (finite) eigenvalues and, without loss of generality,
we may assume that $\,R(\la)\in\partial\mathcal{B}
(P,s_n(\la_0)\,w(|\la_0|)^{-1},w)\subset\textup{Int} [
\mathcal{B}(P,\eps,w)]\,$ (see Lemma \ref{lem.bs}). Then $\la_0$
is also an eigenvalue of all matrix polynomials
\[
  R_{\alpha}(\la)\,=\, ( R_m + \alpha I) \la^m + R_{m-1} \la^{m-1}
                  + \cdots + R_1 \la + R_0 - ( \alpha \la_0^m) I
                   \; ; \;\; \alpha \in \CC \setminus \{ 0 \} ,
\]
where $\,\det (R_m +\alpha I)\ne 0\,$ and $R_{\alpha}(\la)\in
\mathcal{B}(P,\eps,w)$ for sufficiently small $| \alpha |$, i.e.,
$\la_0$ lies in $\Lambda_0$. By Corollary~\ref{t2}, $\lep$ does
not have more connected components than $\,\Lambda_0 \; (\subseteq
\lep)$. Hence, $\lep$ has at most $k$ connected components. \qed

\medskip

In this theorem, recall that since the leading coefficient of
$P(\la)$ is nonsingular, for $\eps$ sufficiently small, $\lep$ has
exactly $k$ bounded connected components. Thus, our upper bound
for the number of connected components of $\lep$ is always
attainable when $\lep$ is bounded.

\medskip

\begin{propos}\label{prop.use}
If $\lep$ is bounded, then any $\,Q(\la)\in\B(P,\eps,w)\,$ has an eigenvalue
in each of these components. Furthermore, $P(\la)$ and $Q(\la)$
have the same number of eigenvalues (counting algebraic
multiplicities) in each connected component of $\lep$.
\end{propos}
\noindent \textbf{Proof.}  See Theorem 2.3 of \cite{LP}. \qed

%%%%%%%%%%%%%%%%%%%%%%%%%%%%%%%%%%%%%%%%%%%%%%%%%%%%%%%%%%%%%%%%%%%%%%
\section{Multiple eigenvalues of perturbations} \label{section:multiple}
%%%%%%%%%%%%%%%%%%%%%%%%%%%%%%%%%%%%%%%%%%%%%%%%%%%%%%%%%%%%%%%%%%%%%%

In this section, we obtain necessary conditions for the existence
of perturbations of $P(\la)$ with multiple eigenvalues. However, we
first construct two perturbations of $P(\la)$ in
$\mathcal{B}(P,\eps,w)$, which are of special interest. They are
used in an argument generalising that of Alam and Bora (Theorem
4.1 of \cite{AB}) for the standard eigenvalue problem.

Suppose that for a $\,\mu\in\lep\setminus\sigma(P),\,$ the
(nonzero) minimum singular value of the matrix $P(\mu)$ has
multiplicity $\,k\geq 1$. Let also
\[
   s_1(\mu) \,\ge\, s_2(\mu) \,\ge\, \cdots \,\ge\, s_{n-k}(\mu)
   \,>\, s_{n-k+1}(\mu) \,=\, \cdots \,=\, s_n(\mu)  > 0
\]
be the singular values of $P(\mu)$ with associated left singular
vectors $\,u_1 , u_2 , \dots , u_n\,$ and associated right
singular vectors $\,v_1 , v_2 , \dots , v_n$. These singular vectors
satisfy the relations $P(\mu) v_j \,=\, s_j (\mu) \, u_j$ for
$j = 1 , 2 , \dots , n$.

Define the $\,n\times n\,$ unitary matrix $\,\hat{Z} = \left [ u_1
\; u_2 \,\cdots\; u_n \right ] \,\left [ v_1 \; v_2 \,\cdots\; v_n
\right ]^*\,$ and the $\,n\times n\,$ matrix $\,\tilde{Z} = \left
[ u_{n-k+1} \; u_{n-k+2} \,\cdots\; u_n \right ] \,\left [
v_{n-k+1} \; v_{n-k+2} \,\cdots\; v_n \right ]^*\,$ of rank $k$.
Then $\,\hat{Z} v_j = u_j\,$ for all $\,j=1,2,\dots ,n$,
$\,\tilde{Z} v_j = u_j\,$ for all $\,j = n-k+1, n-k+2, \dots , n
,\,$ and $\,\tilde{Z} v_j = 0\,$ for all $\,j = 1, 2, \dots ,
n-k$. Furthermore, the (nonsingular) matrix $\,\hat{E} = - \,
s_n(\mu)\, \hat{Z}\,$ satisfies \eq\label{eq:hat1}
   (P(\mu)+\hat{E}) v_j \,=\, s_n(\mu) \, u_j - s_n(\mu) \, u_j
              \,=\, 0 \; ; \;\; j = n-k+1 , n-k+2 , \dots , n
\en
and
\eq\label{eq:hat2}
   (P(\mu)+\hat{E}) v_j \,=\, s_j(\mu) \, u_j
       - s_n(\mu) \, u_j \,\ne\, 0 \; ; \;\;
       j=1,2,\dots,n-k.
\en
Similarly, the (rank $k$) matrix $\,\tilde{E}= - \,s_n(\mu)\,
\tilde{Z}\,$ satisfies \eq\label{eq:tilde1}
 (P(\mu)+\tilde{E}) v_j =\,0\; ;\;\; j=n-k+1,n-k+2,\dots,n
\en
and
\eq\label{eq:tilde2}
   (P(\mu)+\tilde{E}) v_j \,=\, s_j(\mu) \, u_j
    \,\ne\, 0 \; ; \;\; j = 1 , 2 , \dots , n-k .
\en Note also that $\,\|\hat{E}\|=\|\tilde{E}\|=s_n(\mu)$.

Now define (for a given weight function $w(x)$) the matrices
\[
  \hat{\Delta}_j \,=\, \left ( \frac{\overline{\mu}}{| \mu |}\right )^j w_j
  \, w(|\mu |)^{-1} \hat{E} \; ; \;\; j=0,1,\dots ,m
\]
and
\[
 \tilde{\Delta}_j \,=\, \left ( \frac{\overline{\mu}}{| \mu |} \right )^j w_j
 \, w(|\mu |)^{-1} \tilde{E} \; ; \;\; j=0,1,\dots,m,
\]
where we set $\,\overline{\mu} / |\mu| = 0 \,$ when $\,\mu = 0$.
Then
\[
   \sum_{j=0}^{m} \hat{\Delta}_j \,\mu^j \,=\, \left ( \sum_{j=0}^{m}
   w_j\, | \mu |^j \right ) w(|\mu |)^{-1}\,\hat{E}\,=\,\hat{E}
\]
and
\[
   \sum_{j=0}^{m} \tilde{\Delta}_j \,\mu^j \,=\, \left ( \sum_{j=0}^{m}
   w_j\, | \mu |^j \right ) w(|\mu |)^{-1}\,\tilde{E}\,=\,\tilde{E}.
\]
Thus, for the (full rank) perturbation of $P(\la)$
\eq\label{eq.perturbation1}
   \hat{Q}(\la) \,=\, (P_m + \hat{\Delta}_m) \la^m + \cdots
   + (P_1 + \hat{\Delta}_1) \la + P_0 + \hat{\Delta}_0
\en and the (lower rank) perturbation of $P(\la)$
\eq\label{eq.perturbation2}
   \tilde{Q}(\la) \,=\, (P_m + \tilde{\Delta}_m) \la^m + \cdots
   + (P_1 + \tilde{\Delta}_1) \la + P_0 + \tilde{\Delta}_0 ,
\en we have $\,\hat{Q}(\mu) = P(\mu)+\hat{E}\,$ and
$\,\tilde{Q}(\mu) = P(\mu)+\tilde{E}$. From (\ref{eq:hat1}),
(\ref{eq:hat2}), (\ref{eq:tilde1}) and (\ref{eq:tilde2}), it is
clear that $\mu$ is an eigenvalue of the matrix polynomials
$\hat{Q}(\la)$ and $\tilde{Q}(\la)$ with geometric multiplicity
exactly $k$ and associated right eigenvectors $\,v_{n-k+1},
v_{n-k+2},\dots,v_n$.

Moreover, for every $\,j = 0,1,\dots ,m$,
\[
   \|\hat{\Delta}_j\| \,=\, w_j\, w(|\mu |)^{-1}\|\hat{E}\|
   \,=\; \frac{w_j \, s_n(\mu)}{w(|\mu |)}
   \;\le\, \eps\, w_j
\]
and
\[
   \|\tilde{\Delta}_j\| \,=\, w_j\, w(|\mu |)^{-1} \|\tilde{E}\|
   \,=\; \frac{w_j \, s_n(\mu)}{w(|\mu |)}
   \;\le\, \eps\, w_j .
\]
Consequently, $\hat{Q}(\la)$ and $\tilde{Q}(\la)$ lie in
$\mathcal{B}(P,\eps,w)$ and the next result follows:

\begin{propos}  \label{propos.geometric}
Let $\mu \in \lep \setminus \sigma(P)$ and let the nonzero
singular value $\,s_n(\mu)$ $(\le \eps\,w(|\mu |))$ of the matrix
$P(\mu)$ have multiplicity $\,k \geq 1$. Then the perturbation
$\hat{Q}(\la)$ in \textup{(\ref{eq.perturbation1})} and the
perturbation $\tilde{Q}(\la)$ in
\textup{(\ref{eq.perturbation2})} lie in $\mathcal{B}(P,\eps,w)$
and have $\mu$ as an eigenvalue with geometric multiplicity equal
to $k$.
\end{propos}

Clearly, every fault point of $P(\la)$ in $\CC \setminus
\sigma(P)$ is a multiple eigenvalue of $\hat{Q}(\la)$ and
$\tilde{Q}(\la)$ with geometric multiplicity greater than $1$.
Furthermore, in the above discussion, note that for every $\,j =
n-k+1, n-k+2, \dots , n$,
\[
    u^*_j P(\mu) \,=\, s_n(\mu) \, v^*_j
\]
and
\begin{eqnarray*}
  u^*_j(P(\mu)+\hat{E}) &=& u^*_j P(\mu) - s_n(\mu)\, u^*_j \hat{Z} \\
  &=& s_n(\mu)\, v^*_j - s_n(\mu) (\hat{Z}^* u_j)^*    \\
  &=& s_n(\mu)\, v^*_j - s_n(\mu)\, v^*_j \;\,=\;\, 0 .
\end{eqnarray*}
Similarly, for every $\,j = n-k+1, n-k+2 , \dots ,
n ,\,$ we have $\,u^*_j(P(\mu)+\tilde{E})=0$. Thus, $\,u_{n-k+1} ,
u_{n-k+2} , \dots , u_n \,$ are left eigenvectors of the
perturbations $\hat{Q}(\la)$ and $\tilde{Q}(\la)$ in
(\ref{eq.perturbation1}) and (\ref{eq.perturbation2}),
corresponding to $\mu$.

The perturbations $\hat{Q}(\la)$ and $\tilde{Q}(\la)$ defined by
(\ref{eq.perturbation1}) and (\ref{eq.perturbation2}) depend on
$w(x)$ (which is considered fixed) and on the choice of $\mu$. It
is also worth noting that for $\mu=0\,$ and a given weight
function $w(x)$ with a constant coefficient $\,w_0>0,\,$ the
construction of $\hat{Q}(\la)$ and $\tilde{Q}(\la)$ is independent
of the non-constant part of $w(x)$ and requires only $w_0$. In the
remainder of this paper, and without loss of generality, for the
definition of $\hat{Q}(\la)$ and $\tilde{Q}(\la)$, we use the
\textit{constant weight function} $\,w_c(x)=w_0$ $(>0)$ instead of
$w(x)$ whenever $\,\mu=0$.

Using Lemma \ref{lem.con}, one can estimate the (spectral norm)
distance from $P(\la)$ to the set of matrix polynomials that have
a prescribed $\,\mu \notin \sigma(P)\,$ as an eigenvalue (cf.
Lemma 3 of \cite{Tis}).

\begin{cor} \label{cor.eigdistance}
Suppose $\,\mu \notin \sigma(P),\,$ and let $\,\delta = s_n(\mu)
\, w(|\mu|)^{-1}$. Then the perturbations $\hat{Q}(\la)$ and
$\tilde{Q}(\la)$ lie on $\partial \mathcal{B}(P,\delta,w)$ and
have $\mu$ as an eigenvalue. Moreover, for every $\,\eps < \delta
,\,$ no perturbation of $P(\la)$ in $\mathcal{B}(P,\eps,w)$ has
$\mu$ as an eigenvalue.
\end{cor}

\begin{propos}  \label{propos.multi}
Let $\,\mu \in \sigma(P),\,$ and let $\,u,v \in \CC^n\,$ be left
and right eigenvectors of $P(\la)$ corresponding to $\mu$,
respectively. If the derivative of $\,P(\la)$ satisfies $\,u^*
P'(\mu) v = 0,\,$ then $\mu$ is a multiple eigenvalue of $P(\la)$.
\end{propos}

\noindent \textbf{Proof.} If the geometric multiplicity of
$\,\mu\in\sigma(P)\,$ is greater than $1$, then the proposition
obviously holds. Hence, we assume that $\mu$ is an eigenvalue of
$P(\la)$ with geometric multiplicity $1$. For every vector $\,y
\in \CC^n ,\,$ $u^*P(\mu)y = 0,\,$ and thus, $\,u \perp
\textup{Range}[P(\mu)]$. Since $\,u \perp P'(\mu) v\,$ and the
dimension of $\textup{Range}[P(\mu)]$ is $\,n-1,\,$ it follows
that the vector $P'(\mu)v$ belongs to $\textup{Range}[P(\mu)]$,
i.e., there exists a $\,y_{\mu} \in \CC^n \,$ such that
\[
    P(\mu) y_{\mu} + P'(\mu) v \,=\, 0 .
\]
This shows that $\mu$ is a multiple eigenvalue of $P(\la)$ with
the Jordan chain $\{v,\,y_\mu\}$ (see \cite{GLR} for properties of
Jordan chains of matrix polynomials). This implies that $\mu$
is a defective multiple eigenvalue of $P(\la)$. \qed
\medskip

Recall the function $\,F_{\eps}(x,y)\equiv F_{\eps}(x+iy)$ ($x,y
\in\RR$) defined in (\ref{e20}).

\begin{propos}  \label{propos.zerograd}
Suppose that for a point $\,\mu = x_{\mu} + i y_{\mu}\,$ of
$\,\lep \setminus \sigma(P),\,$ $s_n(\mu)$ is a simple singular
value of $P(\mu)$ and $u_{\mu}, v_{\mu}$ are associated left and
right singular vectors, respectively, assuming that
$\,w(x)=w_c(x)\;(=w_0>0)\,$ when $\,\mu=0$. Let $\,\delta =
s_n(\mu)\, w(|\mu|)^{-1}$ $(\le\eps )\,$ and consider the
perturbations $\,\hat{Q}(\la),\,\tilde{Q}(\la) \in
\partial \mathcal{B}(P,\delta,w)\,$ defined by
\textup{(\ref{eq.perturbation1})} and
\textup{(\ref{eq.perturbation2})}. If the
gradient of the function $\, F_{\delta}(x,y) \equiv F_{\delta} (x
+ i y)$ at $\mu$ is zero, then $\mu$ is a defective eigenvalue of
$\hat{Q}(\la)$ and $\tilde{Q}(\la)$ with geometric multiplicity
$1$.
\end{propos}

\noindent \textbf{Proof.} Suppose $\,\mu\ne 0,\,$ and let $\,\nabla
F_{\delta}(x_{\mu},y_{\mu}) = 0,\,$ or equivalently (see Lemma
\ref{lem.an}), let
\[
   \textup{Re} \left ( u_{\mu}^* \, \frac{\partial P(\mu)}
   {\partial x} \, v_{\mu} \right ) \,=\,
   \delta \, \frac{\partial w(|\mu|)}{\partial x}
   \;\; \mbox{ and } \;\;
   \textup{Re} \left ( u_{\mu}^* \, \frac{\partial P(\mu)}
   {\partial y} \, v_{\mu} \right ) \,=\,
   \delta \, \frac{\partial w(|\mu|)}{\partial y} \, .
\]
Since
\[
  \frac{\partial P(\mu)}{\partial x} \;=\, P'(\mu)
  \;\; \mbox{ and } \;\;
  \frac{\partial P(\mu)}{\partial y} \;=\, i \, P'(\mu) ,
\]
we see that
\[
  \textup{Im} \left ( u_{\mu}^* \, \frac{\partial P(\mu)}
  {\partial x} \, v_{\mu} \right ) \,=\, - \,
  \textup{Re} \left ( u_{\mu}^* \, \frac{\partial P(\mu)}
  {\partial y} \, v_{\mu} \right ) .
\]
Moreover,
\[
   \frac{\partial w(|\mu|)}{\partial x} \;=\;
   \frac{x_{\mu}}{|\mu|} \; w'(|\mu|)
   \;\; \mbox{ and } \;\;
   \frac{\partial w(|\mu|)}{\partial y} \;=\;
   \frac{y_{\mu}}{|\mu|} \; w'(|\mu|) ,
\]
and consequently,
\[
   u_{\mu}^* P'(\mu) v_{\mu} \,=\,
   u_{\mu}^* \, \frac{\partial P(\mu)}
   {\partial x} \, v_{\mu} \,=\,
   \delta \,\frac{\partial w(|\mu|)}{\partial x} - \,i\,
   \delta \,\frac{\partial w(|\mu|)}{\partial y} \,=\,
   \delta \, \frac{\overline{\mu}}{|\mu|} \; w'(|\mu|) .
\]

Consider the perturbation
\[
   \hat{Q}(\la) \,=\, (P_m + \hat{\Delta}_m) \la^m + \cdots
   + (P_1 + \hat{\Delta}_1) \la + P_0 + \hat{\Delta}_0
\]
in (\ref{eq.perturbation1}). Then $\hat{Q}(\la)$ lies on the
boundary of the (compact) set $\,\mathcal{B}(P,\delta,w) \subseteq
\mathcal{B}(P,\eps,w)\,$ and its derivative satisfies
\begin{eqnarray*}
   u_{\mu}^* \hat{Q}'(\mu) v_{\mu} &=& u_{\mu}^*
   P'(\mu) v_{\mu} + u_{\mu}^* \left ( \sum_{j=1}^{m} j \,
   \hat\Delta_j \, \mu^{j-1} \right )  v_{\mu} \\
   &=& \delta \,\frac{\overline{\mu}}{|\mu|}\; w'(|\mu|)
   + ( u_{\mu}^* \hat{E} v_{\mu}) \; \frac{w'(|\mu|)}
   {w(|\mu|)} \; \frac{\overline{\mu}}{|\mu|}  \\
   &=& \delta \,\frac{\overline{\mu}}{|\mu|}\; w'(|\mu|)
   - \, \frac{s_n(\mu)}{w(|\mu|)} \,
   \frac{\overline{\mu}}{|\mu|} \; w'(|\mu|) \\
   &=& \delta \, \frac{\overline{\mu}}{|\mu|} \; w'(|\mu|)
   - \delta \, \frac{\overline{\mu}}{|\mu|} \; w'(|\mu|)
   \;\, = \;\, 0 ,
\end{eqnarray*}
where $u_{\mu}$ and $v_{\mu}$ are left and right eigenvectors of
$\hat{Q}(\la)$ corresponding to $\mu$, respectively (see
Proposition \ref{propos.geometric} and the related discussion).
The same is also true for the perturbation $\tilde{Q}(\la)$ in
(\ref{eq.perturbation2}) and its derivative. By Propositions
\ref{propos.geometric} and \ref{propos.multi}, $\mu$ is a multiple
eigenvalue of $\hat{Q}(\la)$ and $\tilde{Q}(\la)$ with geometric
multiplicity $1$.

For $\,\mu =0,\,$ the proof is the same, keeping in mind that the
constant weight function $\,w_c(x)=w_0\;(>0)\,$ is differentiable
(with zero partial derivatives) at the origin. \qed

%%%%%%%%%%%%%%%%%%%%%%%%%%%%%%%%%%%%%%%%%%%%%%%%%%%%%%%%%%%%%
\section{Multiple points on $\partial\lep$ and connected components
of $\lep$} \label{sec.7}
%%%%%%%%%%%%%%%%%%%%%%%%%%%%%%%%%%%%%%%%%%%%%%%%%%%%%%%%%%%%%

At first glance it may seem that multiple (crossing) points on
$\partial \lep $  will be exceptional. However, when we consider the
evolution of $\partial \lep$ as $\eps$ increases, it is clear that,
as disjoint components of $\lep$ expand, there will be critical
values of $\eps$ at which they meet and multiple points are created.

Next, based on the results of the previous section, we show that
multiple points of $\partial\lep$ are multiple eigenvalues of
perturbations of $P(\la)$ on $\partial \mathcal{B}(P,\eps,w)$ and,
also, these perturbations can be constructed explicitly. (Recall
that, when $\mu =0$, we use the constant weight function
$\,w_c(x)=w_0>0\,$ for the definition of the perturbations
$\hat{Q}(\la)$ and $\tilde{Q}(\la)$ in (\ref{eq.perturbation1})
and (\ref{eq.perturbation2}).)

\begin{thm} \label{thm.inter}
Suppose that, as the parameter $\,\eps > 0\,$ increases, two
different connected components of $\,\lep\ne\CC ,\,$ $\GG_1$ and
$\GG_2$, meet at $\,\mu\in\CC$. Then the following hold:

\begin{description}
\item[(i)] If $\,\mu\ne 0,\,$ then it is a multiple eigenvalue of
           the perturbations $\,\hat{Q}(\la),\tilde{Q}(\la)\in\partial
           \mathcal{B}(P,\eps,w)\,$ defined by
           \textup{(\ref{eq.perturbation1})} and
           \textup{(\ref{eq.perturbation2})}.
\item[(ii)] If $\,\mu =0\,$ and $\,w(x)=w_c(x)\; (=w_0>0),\,$ then
            $\,\mu=0\,$ is a multiple eigenvalue of the perturbations
            $\,\hat{Q}(\la), \tilde{Q}(\la) \in \partial
            \mathcal{B}(P,\eps,w_c)$.
\item[(iii)] If $\,\mu =0$, $\,w(x)\ne w_c(x),\,$ $\lep$ is
             bounded and the origin is the only intersection point of
             $\,\GG_1$ and $\GG_2$, then $\,\mu=0\,$ is a multiple eigenvalue
             of a perturbation on $\partial\mathcal{B}(P,\eps,w)$.
\end{description}

Furthermore, in the first two cases, if $s_n(\mu)$ is a simple
singular value of $P(\mu)$, then $\mu$ is a defective eigenvalue
of $\hat{Q}(\la)$ and $\tilde{Q}(\la)$ with geometric
multiplicity $1$.
\end{thm}

\noindent \textbf{Proof.} Suppose that $\,s_n(\mu)$ $(=\eps\,
w(|\mu|))\,$ is a multiple singular value of the matrix $P(\mu)$.
Then by Proposition \ref{propos.geometric}, the perturbations
$\,\hat{Q}(\la), \tilde{Q}(\la) \in \partial
\mathcal{B}(P,\eps,w)\,$ have $\mu$ as a multiple eigenvalue of
geometric multiplicity greater than $1$. Hence, we may assume that
$s_n(\mu)$ is a simple singular value of $P(\mu)$, and consider
the three cases of the theorem.

\smallskip

\noindent (i) Suppose $\,\mu \ne 0,\,$ and recall (\ref{eq.bo}).
By virtue of Lemma \ref{lem.an}, $F_\varepsilon(x,y)$ is real
analytic in a neighbourhood of $\mu$. Furthermore, $\nabla
F_\varepsilon(\mu)=0$, otherwise the implicit function theorem
would ensure the existence of a smooth curve on a neighbourhood of
$\mu$ parameterising $\partial \lep$ and contradict the
fact that $\partial \GG_1 \cap
\partial \GG_2$ is a finite set (Theorem~\ref{p7}).
Therefore, Proposition \ref{propos.zerograd} yields the desired
conclusion.

\smallskip

\noindent (ii) If $\,\mu=0\,$ and $\,w(x)=w_c(x)\;(=w_0>0),\,$
then the result follows by applying Proposition
\ref{propos.zerograd} as in case (i).

\smallskip

\noindent (iii) Suppose $\lep$ is bounded, $\,w(x)\ne w_c(x),\,$
and $\,\mu =0\,$ is the only intersection point of $\GG_1$ and
$\GG_2$. By Proposition \ref{prop.use}, for any positive $\,\delta
< \eps,\,$ all the perturbations in $\mathcal{B}(P,\delta,w)$ have
a constant number of eigenvalues in
$\,\Lambda_{\delta}(P)\cap\GG_j ,\,$ say $\kappa_j$, for
$\,j=1,2$. Here and throughout this proof, eigenvalues are
counted according to their algebraic multiplicities.

Define the sets
\[
  \mathcal{B} \,=\, \left\{ Q(\la) \in \mathcal{B}(P,\eps,w)  :\,
  0 \in\sigma(Q) \right\} \,\subseteq\, \partial \mathcal{B}(P,\eps,w)
\]
and
\[
   \mathcal{B}_j \,=\, \left \{ Q(\la)\in\mathcal{B}:\,Q(\la)
   \;\,\mbox{has less than}\;\,\kappa_j \;\,\mbox{eigenvalues in}
   \;\,\GG_j \setminus \{ 0 \} \right \} \, ; \;\; j=1,2.
\]
If $\,\mathcal{B}_j = \emptyset\,$ ($j=1,2$), then Proposition
\ref{prop.use} and the continuity of the eigenvalues of matrix
polynomials with respect to the entries of their coefficients
imply that $\,0 \notin \GG_j$; this is a contradiction. Hence, the
sets $\mathcal{B}_1$ and $\mathcal{B}_2$ are both non-empty.

Now consider the constant weight function $\,w_c(\la) = w_0 \;
(>0)\,$ and the associated $\,\eps$-pseudospectrum of $P(\la)$,
\[
   \Lambda_{\eps,w_c}(P) \,=\, \left\{\mu\in\CC :\,
   \det Q(\mu) = 0 ,\; \| \Delta_0 \| \le \eps\, w_0 ,\;
   \Delta_1 = \cdots = \Delta_m = 0 \right \} .
\]
Clearly, $\,\Lambda_{\eps,w_c}(P)\subseteq\lep\,$ and $\,0 \in
\partial \Lambda_{\eps,w_c}(P)$. For any $\,j=1,2,\,$ consider a
perturbation
\[
  Q_j(\la)\,=\,(P_m + \Delta_m) \la^m + \cdots +
               (P_1 + \Delta_1) \la + P_0 + \Delta_0
               \; \in \; \mathcal{B}_j ,
\]
and define the matrix polynomial
\[
    Q_{j,c}(\la) \,=\, P_m \la^m + \cdots + P_1 \la + P_0
                    + \Delta_0 \,=\, P(\la) +  \Delta_0
    \;\in\; \mathcal{B} \cap \partial \mathcal{B}(P,\eps,w_c)
\]
and the continuous trajectory
\[
     Q_j(t;\la)\,=\,(P_m + t\Delta_m) \la^m + \cdots +
                (P_1 + t\Delta_1) \la + P_0 + \Delta_0
                \; \in \; \mathcal{B} \; ; \;\; 0 \le t \le 1
\]
with $\,Q_j(0;\la)=Q_{j,c}(\la)\,$ and $\,Q_j(1;\la)=Q_j(\la)$. If
$\,\mu=0\,$ is a multiple eigenvalue of $Q_j(t;\la)$ for some
$\,t\in[0,1],\,$ then there is nothing to prove.

Let $\,\mu=0\,$ be a simple eigenvalue of $\,Q_j(t;\la)\in
\mathcal{B}\,$ for all $\,t\in[0,1]$. Since $\lep$ is bounded and
the origin is the only intersection point of $\GG_1$ and $\GG_2$,
by the continuity of the eigenvalues with respect to the
coefficient matrices, it follows that all $Q_j(t;\la)$ ($0\le t
\le 1$) have exactly $\,\kappa_j-1\,$ eigenvalues in $\,\GG_j
\setminus\{0\},\,$ i.e., they lie in $\mathcal{B}_j$. Thus,
$\,Q_{j,c}(\la)\in\mathcal{B}_j$. Moreover, again by Proposition
\ref{prop.use} and the continuity of eigenvalues, an eigenvalue of
the matrix polynomials $\,P(\la)+(1-t)\Delta_0\,$ ($0\le t\le 1$)
traces a continuous path in $\GG_j$ connecting the origin with an
eigenvalue of $P(\la)$. This means that the origin is an
intersection point of $\Lambda_{\eps,w_c}(P)\cap\GG_1$ and
$\Lambda_{\eps,w_c}(P)\cap\GG_2$. Hence, $\,\mu=0\,$ is a
multiple point of $\partial\Lambda_{\eps,c}(P)$, and as in (ii),
it is a multiple eigenvalue of the perturbations $\,\hat{Q}(\la),
\tilde{Q}(\la)\in\partial\mathcal{B}(P,\eps,w_c)\subset
\partial\mathcal{B}(P,\eps,w)$. \qed

\medskip

Now we can generalise a theorem of Mosier concerning scalar
polynomials (Theorem 3 of \cite{Mo}).

\begin{thm} \label{thm.multi}
Suppose $\lep$ is bounded and $\,\GG$ is a connected component of
$\lep$. Then the matrix polynomial $P(\la)$ has more than one
eigenvalue in $\,\GG$ (counting multiplicities) if and only if
there is a perturbation $\,Q(\la) \in \mathcal{B}(P,\eps,w)\,$
with a multiple eigenvalue in $\,\GG$.
\end{thm}

\noindent \textbf{Proof.} For the converse part, it is clear that
if a perturbation $\,Q(\la) \in \mathcal{B}(P,\eps,w)\,$ has a
multiple eigenvalue in $\GG$, then by Proposition \ref{prop.use},
$P(\la)$ has at least two eigenvalues in $\GG$, counting
multiplicities.

For the sufficiency, if the matrix polynomial $P(\la)$ has a
multiple eigenvalue in $\GG$, then there is nothing to prove.
Thus, we assume that $P(\la)$ has two simple eigenvalues, $\la_1$
and $\la_2$, in $\GG$. By the continuity of the eigenvalues with
respect to the coefficient matrices, it follows that there is a
positive $\,\delta \le \eps ,\,$ such that $\Lambda_{\delta}(P)$
has a (bounded) connected component $\,\GG_{\delta} \subseteq
\GG\,$ that is composed of two compact sets, $\GG_{1,\delta}$ and
$\GG_{2,\delta}$, with disjoint interiors and intersecting
boundaries. Moreover, without loss of generality, we can assume
that $\la_1$ and $\la_2$ lie in the interior of $\GG_{1,\delta}$
and $\GG_{2,\delta}$, respectively. Then the curve enclosing
$\GG_{\delta}$ either crosses itself or is tangent to itself at
some point $\,\mu_{\delta}\in \CC \setminus \sigma(P)$. The result
follows from Theorem~\ref{thm.inter}. Note that if $\GG_{1,\delta}$ and
$\GG_{2,\delta}$ intersect at the origin and one other point, 
Theorem~\ref{thm.inter} does not apply at $\mu=0$, but it will at
that other point. \qed

%%%%%%%%%%%%%%%%%%%%%%%%%%%%%%%%%%%%%%%%%%%%%%%%%%%%%%
\section{Two numerical examples} \label{examples}
%%%%%%%%%%%%%%%%%%%%%%%%%%%%%%%%%%%%%%%%%%%%%%%%%%%%%%

We present two numerical examples, which illustrate the results of
the previous section and suggest possible applications. The
figures were drawn using the boundary-tracing algorithm described
in \cite{LP}.

\medskip

\begin{figure}
\begin{center}
   \epsfig{file=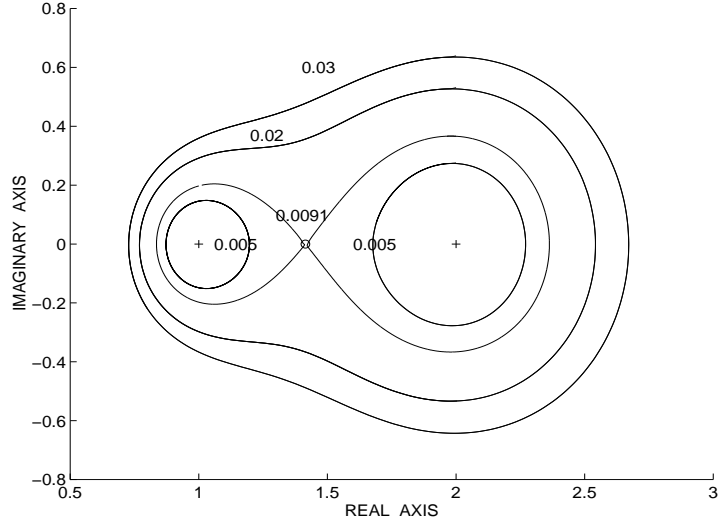, height=70mm, width=95mm, clip=}
   \caption{\small A single intersection point.}
   \label{figure1}
\end{center}
\end{figure}

\noindent \textbf{Example 6} The spectrum of the $\,2 \times 2\,$
quadratic matrix polynomial
\[
     P (\la) \,=\,  \left [ \begin{array}{cc}
                          (\la - 1)^2 &  \la  \\
                                0     & (\la - 2)^2
                              \end{array} \right ]
             \,=\,  \left [ \begin{array}{cc}
                                1 & 0  \\
                                0 & 1
                              \end{array} \right ] \la^2 +
        \left [ \begin{array}{cc}
                                -2 &  1  \\
                                 0 & -4
                              \end{array} \right ] \la +
        \left [ \begin{array}{cc}
                                1 &  0  \\
                                0 &  4
                               \end{array} \right ]
\]
is $\,\sigma(P) = \{ 1 , 2 \}$. Both eigenvalues are plotted as
``+'' in Figure \ref{figure1}, and have algebraic multiplicity
equal to $2$ and geometric multiplicity equal to $1$. The
boundaries $\partial \lep$ for $\,w(x) = x^2 + x + 1\,$ (i.e., for
perturbations measured in the absolute sense) and for $\,\eps =
0.005 ,\, 0.0091 ,\, 0.02 ,\, 0.03,\,$ are also sketched in
Figure \ref{figure1}.

Assuming that the pseudospectrum $\Lambda_{0.0091}(P)$ is
connected with one self-intersection $\,\mu = 1.4145\,$ plotted as ``o'',
this figure indicates that $\lep$ consists of two connected
components for $\,\eps < 0.0091 ,\,$ and that it is connected for
$\,\eps \geq 0.0091$. Moreover, the singular values of the matrix
$P(\mu)$ are $s_1(\mu) = 1.4650$ and $s_2(\mu) = 0.0402 ,\,$
i.e., $s_2(\mu)$ is simple ($\mu$ is not a fault point of
$P(\la)$) and the function
\[
         F_{0.0091}(x,y) \,\equiv\, F_{0.0091}(x + i y)
                     \,=\, s_2(x + i y) - 0.0091
               \, w( | x + i y | ) \; ; \;\; x,y \in \RR
\]
has zero gradient at the point $\mu$. Thus, by Proposition
\ref{propos.zerograd} and Theorem \ref{thm.inter}, two
perturbations of $P(\la)$ on the boundary of $\mathcal{B}
(P,0.0091,w)$ that have $\mu$ as a defective eigenvalue are
$\hat{Q}(\la)$ and $\tilde{Q}(\la)$ in (\ref{eq.perturbation1})
and (\ref{eq.perturbation2}), and can be easily constructed.
%%%%%%%%%%%%%%%%%%%%%%%%%%%%%%%%%%%%%%%%%%%%%%%%%%%
%%%%%%%%%%%%%%%%%%%%%%%%%%%%%%%%%%%%%%%%%%%%%%%%%%%
Left and right singular vectors of $P(\mu)$ corresponding to
$s_1(\mu)$ are
\[
   u_1 \,=\, \left [ \begin{array}{c} 0.9726 \\
                         0.2325 \end{array} \right ]
     \;\; \mbox{ and } \;\;
   v_1 \,=\, \left [ \begin{array}{c} 0.1141 \\
                         0.9935 \end{array} \right ] ,
\]
respectively, and left and right singular vectors of $P(\mu)$
corresponding to $s_2(\mu)$ are
\[
   u_2 \,=\, \left [ \begin{array}{c} -0.2325 \\
                         0.9726 \end{array} \right ]
     \;\; \mbox{ and } \;\;
   v_2 \,=\, \left [ \begin{array}{c} -0.9935 \\
                         0.1141 \end{array} \right ] ,
\]
respectively.

The unitary matrix
\[
     \hat{Z} \,=\, \left [ \begin{array}{cccc}
                   u_1 & u_2 \end{array} \right ]
                   \left [ \begin{array}{cccc}
                   v_1 & v_2 \end{array} \right ]^*
     \,=\, \left [ \begin{array}{cc} 0.3419 & 0.9397 \\
                       -0.9397 & 0.3419 \end{array} \right ]
\]
satisfies
\[
   \hat{Z} v_1 = u_1 , \;\; u_1^* \hat{Z} = v_1^* , \;\;
   \hat{Z} v_2 = u_2 \;\; \mbox{ and } \;\; u_2^* \hat{Z} = v_2^* ,
\]
and the rank one matrix
\[
     \tilde{Z} \,=\, u_2 v_2^*
     \,=\, \left [ \begin{array}{cc} 0.2310 & -0.0265 \\
                       -0.9663 & 0.1110 \end{array} \right ]
\]
satisfies
\[
   \tilde{Z} v_1 = 0 , \;\;  u_1^* \tilde{Z} = 0 , \;\;
   \tilde{Z} v_2 = u_2 \;\; \mbox{ and } \;\; u_2^* \tilde{Z} = v_2^* .
\]

We define the matrices
\[
   \hat{\Delta}_0 \,=\, \hat{\Delta}_1 \,=\, \hat{\Delta}_2
   \,=\, (\mu^2 + \mu +1)^{-1} ( - s_2(\mu) \hat{Z} )
   \,=\, \left [ \begin{array}{cc} - 0.0031 & - 0.0086 \\
          0.0086 & - 0.0031 \end{array} \right ]
\]
and the matrices
\[
   \tilde{\Delta}_0\,=\,\tilde{\Delta}_1\,=\,\tilde{\Delta}_2
   \,=\, (\mu^2 + \mu +1)^{-1} ( - s_2(\mu) \tilde{Z})
   \,=\, \left [ \begin{array}{cc} - 0.0021 & 0.0002 \\
          0.0088 & - 0.0010 \end{array} \right ] ,
\]
all with spectral norm $\,0.0091$. Then the perturbations
\[
  \hat{Q} (\la) \,=\, P(\la) + (\hat{\Delta}_2 \la^2
  + \hat{\Delta}_1 \la + \hat{\Delta}_0)
\]
\[
  \,=\, \left [ \begin{array}{cc} 0.9969 & -0.0086 \\
          0.0086 & 0.9969 \end{array} \right ] \la^2
 + \left [ \begin{array}{cc} -2.0031 & 0.9914 \\
         0.0086 & -4.0031 \end{array} \right ] \la
 + \left [ \begin{array}{cc} 0.9969 & -0.0086 \\
                  0.0086 & 3.9969 \end{array} \right ]
\]
and
\[
  \tilde{Q} (\la) \,=\, P(\la) + (\tilde{\Delta}_2 \la^2
  + \tilde{\Delta}_1 \la + \tilde{\Delta}_0)
\]
\[
  \,=\, \left [ \begin{array}{cc} 0.9979 & 0.0002 \\
          0.0088 & 0.9990 \end{array} \right ] \la^2
 + \left [ \begin{array}{cc} -2.0021 & 1.0002 \\
         0.0088 & -4.0010 \end{array} \right ] \la
 + \left [ \begin{array}{cc} 0.9979 & 0.0002 \\
                  0.0088 & 3.9990 \end{array} \right ]
\]
lie on $\partial \mathcal{B}(P,0.0091,w)$ and have a multiple
eigenvalue (approximately) equal to $\,\mu = 1.4145\,$ with
algebraic multiplicity $2$ and geometric multiplicity $1$,
confirming our results. \qed

\medskip

It is important to note that, by Theorems \ref{thm.inter} and
\ref{thm.multi}, pseudospectra yield a visual approximation of the
distance to multiple eigenvalues, i.e., the spectral norm distance
from an $\,n\times n\,$ matrix polynomial $P(\la)$ with a
nonsingular leading coefficient and all its eigenvalues simple to
$n\times n$ matrix polynomials with multiple eigenvalues. For a given weight
function $w(x)$, this distance is defined by
\begin{eqnarray*}
   \textup{r} (P) &:=& \min \{ \eps > 0 : \exists \;\,
                  Q(\la) \in \mathcal{B}(P,\eps,w) \;\,
                  \mbox{with multiple eigenvalues} \}     \\
      &\equiv&   \min \{ \eps > 0 : \exists  \;\,
            Q(\la) \in \partial \mathcal{B}(P,\eps,w) \;\,
            \mbox{with multiple eigenvalues} \} .
\end{eqnarray*}
Then Theorems \ref{thm.inter} and \ref{thm.multi} imply the
following result (see \cite{AB,Mal} for the standard eigenvalue
problem).

\begin{cor} \label{distance}
Let $P(\la)$ be an $\,n\times n\,$ matrix polynomial as in
\textup{(\ref{eq.poly})} with a nonsingular leading coefficient
and simple eigenvalues only.
\begin{description}
\item[(a)] If $\lep$ is bounded, then
\[
   \textup{r}(P) \,=\, \min \{ \eps > 0 : \, \lep \;
  \mbox{has less than} \;\,n m\; \mbox{connected components} \} .
\]
\item[(b)] If $\lep$ is unbounded and, as $\eps$ increases from
zero, its connected components meet at points different from the
origin, then
\[
  \hspace{-0.80cm}
  \textup{r}(P) \,=\, \min \{ \eps > 0 :\, \mbox{the number
  of connected components of} \;\,\lep\; \mbox{decreases} \} .
\]
\end{description}
\end{cor}

\begin{figure}
\begin{center}
   \epsfig{file=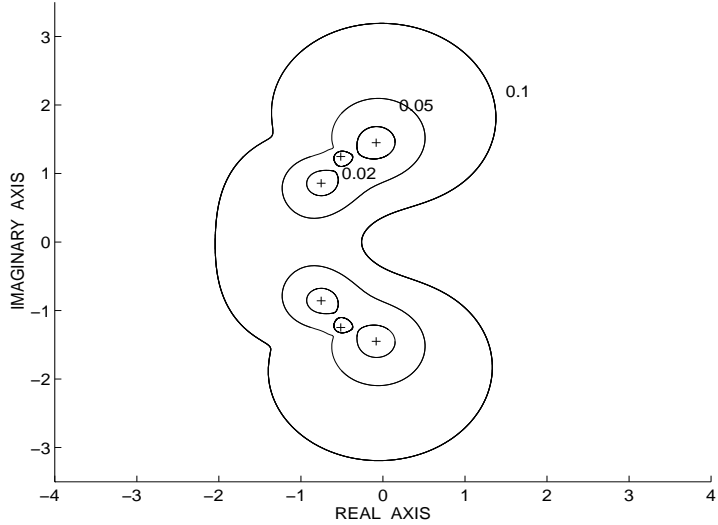, height=70mm, width=95mm, clip=}
   \caption{\small A damped vibrating system.}
   \label{figure2}
\end{center}
\end{figure}

\noindent \textbf{Example 7} Consider the $\,3 \times 3\,$
self-adjoint matrix polynomial
\[
  P(\la) \,=\, \left [ \begin{array}{ccc}
                                1 & 0 & 0  \\
                                0 & 2 & 0  \\
                                0 & 0 & 5
                              \end{array} \right ] \la^2 +
        \left [ \begin{array}{ccc}
                                0 &  0 &  0  \\
                                0 &  3 & -1  \\
                                0 & -1 &  6
                              \end{array} \right ] \la +
        \left [ \begin{array}{ccc}
                                2 & -1 & 0  \\
                               -1 &  3 & 0  \\
                                0 &  0 & 10
                              \end{array} \right ]
\]
(see \cite[Example 5.2]{LP}), which corresponds to a damped
vibrating system. The boundaries of $\lep$ for $\,w(x)=\| A_ 2 \|
x^2 + \| A_ 1 \| x + \| A_ 0 \| = 5 x^2 + 6.3 x + 10\,$ (i.e., for
perturbations measured in a relative sense) and for $\,\eps = 0.02
,\, 0.05 ,\, 0.1 ,\,$ are drawn in Figure \ref{figure2}. The
eigenvalues of $P(\la)$, $\, -0.08 \pm i 1.45, \, -0.75 \pm i 0.86
\,$ and $\, -0.51 \pm i 1.25 , \,$ are plotted as ``+''.
%%%%%%%%%%%%%%%%%%%%%%%%%%%%%%%%%%%%%%%%%%%%%%%%%%%
%%%%%%%%%%%%%%%%%%%%%%%%%%%%%%%%%%%%%%%%%%%%%%%%%%%

We learn from this figure and the above discussion that there
exist an $\,\eps_1 = \textup{r}(P)\,$ in $\,( 0.02 , 0.05 )\,$
(for which, the pseudospectrum starts having less than six
connected components) and an $\eps_2$ in $\,( 0.05 , 0.1 )\,$ (for
which, the pseudospectrum becomes connected) such that the
following hold:
\begin{enumerate}
\item For every $\,\eps < \eps_1 ,\,$ all the perturbations
      $\,Q(\la) \in \mathcal{B}(P,\eps,w)\,$ have
      only simple eigenvalues.
\item For every $\,\eps \in [ \eps_1 , \eps_2 ) ,\,$
      some perturbations $\,Q(\la) \in \mathcal{B}(P,\eps,w)\,$
      have multiple non-real eigenvalues (in a neighbourhood between
      the eigenvalues of $P(\la)$ in the open upper half-plane and
      in a neighbourhood between the eigenvalues of $P(\la)$ in the open
      lower half-plane), but no perturbation in $\mathcal{B}(P,\eps,w)$
      has multiple real eigenvalues.
\item For every $\,\eps \geq  \eps_2 ,\,$ some perturbations
      $\,Q(\la)\in \mathcal{B}(P,\eps,w)\,$
      have multiple real eigenvalues in the interval
      $\, [ -2.1 , -0.2 ]$. \qed
\end{enumerate}

%%%%%%%%%%%%%%%%%%%%%%%%%%%%%%%%%%%%%%%%%%%%%
\section*{Acknowledgements}
%%%%%%%%%%%%%%%%%%%%%%%%%%%%%%%%%%%%%%%%%%%%%%
The work of Peter Lancaster was supported in part by a grant from
the Natural Sciences and Engineering Research Council of Canada.
The work of Lyonell Boulton was supported in part by the Leverhulme Trust.

%%%%%%%%%%%%%%%%%%%%%%%%%%%%%%%%%%%%%%%%%%%%%%%%

% --------------------------------------------------------
\end{document}